\documentclass[10pt,onesite,draft]{article}
\usepackage{amsmath,amsxtra,amssymb,latexsym, amsfonts, indentfirst, color}
\usepackage[mathscr]{eucal}
\newfam\cyrfam

\newenvironment{proof}{%
\par\addvspace{6pt plus3pt minus2pt}%
\noindent{\bfseries\itshape\textit Proof.}\ignorespaces} {%
\if@halmos\halmos\fi
\par\addvspace{6pt plus3pt minus2pt} }
\begin{document}
\parskip 4pt
\large
\setlength{\baselineskip}{15 truept}
\setlength{\oddsidemargin} {0.5in}
\overfullrule=0mm
\def\bfh{\vhtimeb}
\date{}

\title{\bf \large   ON  SOME  MULTIPLICITY AND MIXED
\\MULTIPLICITY FORMULAS $^1$\\
     (Forum Math. 26(2014), 413-442)} 
\def\b{\vntime}
\author{
 Duong Quoc Viet  and Truong Thi Hong Thanh \\
\small Department of Mathematics, Hanoi National University of Education\\
\small 136 Xuan Thuy Street, Hanoi, Vietnam\\
\small Emails: duongquocviet@fmail.vnn.vn \;and\; thanhtth@hnue.edu.vn\\}
  \date{}
\maketitle
\centerline{
\parbox[c]{10.4 cm}{
\small  ABSTRACT: This paper gives the additivity and reduction formulas for  mixed multiplicities of multi-graded modules $M$ and mixed multiplicities of arbitrary ideals, and establishes the recursion formulas for the sum of all the  mixed  multiplicities of $M.$ As an application of these formulas  we get  the recursion formulas for the  multiplicity of  multi-graded Rees modules.
  }}

\vspace{24pt}
\centerline{\Large\bf 1. Introduction}$\\$
Let $(A,\frak{m})$ be  an artinian local ring with  maximal ideal $\frak{m}$, infinite residue field $k = A/\frak{m}.$ Let $S=\bigoplus_{n_1,\ldots,n_d\ge 0}S_{(n_1,\ldots,n_d)}$ $(d > 0)$ be  a finitely generated standard $\mathbb{N}^d$-graded algebra over $A$ (i.e., $S$ is generated over $A$ by elements of total degree 1) and let $M=\bigoplus_{n_1,\ldots,n_d\ge 0}M_{(n_1,\ldots,n_d)}$ be  a finitely generated $\mathbb{N}^d$-graded $S$-module such that $M_{(n_1,\ldots,n_d)}=S_{(n_1,\ldots,n_d)}M_{(0,\ldots,0)}$  for all $n_1,\ldots,n_d \ge 0.$
\footnotetext[1]{\begin{itemize}
 \item[ ] This research was in part supported by a grant from  NAFOSTED.
\item[ ]{\it Mathematics  Subject  Classification} (2010): Primary 13H15. Secondary 13A02, 13A15, 13A30, 14C17.
\item[ ]$ Key\; words \;  and \; phrases:$ Noetherian  ring, mixed multiplicity,  multi-graded  module, filter-regular  sequence.
\end{itemize}}
Throughout this paper, put $S_i=S_{(0,\ldots,{\underbrace{1}_i},\ldots,0)}$
\; for\; all\;  $i=1,\ldots,d$ \;and

      $$\begin{array}{lll}
&S^\triangle&=\bigoplus_{n\ge 0}\;S_{(n,\ldots,n)},\;S_{++}=\bigoplus_{\;n_1,\ldots,n_d> 0}S_{(n_1,\ldots,n_d)},\\
&S^\triangle_+&=\bigoplus_{n> 0}\;S_{(n,\ldots,n)},\;S_+=\bigoplus _{n_1+\cdots+n_d> 0}S_{(n_1,\ldots,n_d)},\\&M^\triangle &=\bigoplus_{n\ge 0}M_{(n,\ldots,n)},\;
{\frak a}:{\frak b}^\infty=\bigcup_{n\ge 0}(\frak a:{\frak b}^n)
.\end{array}$$

\noindent Denote by $\text{Proj}\; S$ the set of the homogeneous prime ideals of $S$ which do not contain $S_{++}$.
Set $\dim M^\triangle = \ell \;\;\text{ and }$\;\;$$\text{Supp}_{++}M=\Big\{P\in \text{Proj}\; S\;|\;M_P\ne 0\Big\}.$$
   By \cite[Remark 3.1]{VM}, $\dim\text{Supp}_{++}M = \ell-1.$ And by \cite[Theorem 4.1]{HHRT},  $\ell_A[M_{(n_1,\ldots,n_d)}]$ is a polynomial of degree $\ell-1$ for all large $n_1,\ldots,n_d.$
  The terms of total degree $\ell-1$ in this polynomial have the form
$$ \sum_{k_1\:+\:\cdots\:+\:k_d\;=\;\ell-1}e(M;k_1,\ldots,k_d)\dfrac{n_1^{k_1}\cdots n_d^{k_d}}{k_1!\cdots k_d!}.$$ Then $e(M;k_1,\ldots,k_d)$ are   called the {\it  mixed multiplicity of type $(k_1,\ldots,k_d)$ of $M$} \cite{HHRT}.
In the case that $(R, \frak n)$  is  a  noetherian   local ring with  maximal ideal $\frak{n},$ $J$ is an $\frak n$-primary ideal, $I_1,\ldots, I_d$ are ideals of $R,$ $N$ is   a finitely generated  $R$-module,  then it is easily seen that
$$ F_J(J,I_1,\ldots,I_d;N) =\bigoplus_{n_0, n_1,\ldots,n_d\ge 0}\dfrac{J^{n_0}I_1^{n_1}\cdots I_d^{n_d}{N}}{J^{n_0+1}I_1^{n_1}\cdots I_d^{n_d}{N}}$$
 is a finitely generated graded $F_J(J,I_1,\ldots,I_d;R)$-module.
Mixed multiplicities of $F_J(J,I_1,\ldots,I_d;N)$ are denoted by $e\big(J^{[k_0+1]},I_1^{[k_1]},\ldots,I_d^{[k_d]};N\big)$ and which are called {\it
mixed multiplicities of  $N$ with respect to ideals $J,I_1,\ldots,I_d$}  (see \cite{MV, Ve}).

Although the problems of expressing the multiplicity of graded modules in terms of mixed multiplicities and  the relationship between  mixed multiplicities  and  Hilbert-Samuel multiplicity  have attracted much attention in past years (the citations will be mentioned in the next sections), the properties similar to that of the Hilbert-Samuel multiplicity (for  instance: the additive property on exact sequences as in \cite[Lemma 17.4.4]{HS} and the additivity and reduction formula \cite[Theorem 17.4.8]{HS} for  mixed multiplicities of $\frak n$-primary ideals...) for mixed multiplicities of arbitrary ideals and multi-graded modules, and other properties, are not yet known.

In the present paper, by a new approach we give additivity and reduction formulas for  mixed multiplicities of multi-graded modules  and mixed multiplicities of arbitrary ideals. And  we establish the recursion formulas for the sum of all the  mixed multiplicities of multi-graded modules.

 As one might expect,  we first obtain the following result for mixed multiplicities of  multi-graded modules.

\vskip 0.2cm

\noindent
{\bf Theorem 3.1.}\;{\it Let $S$  be  a finitely generated standard \;$\mathbb{N}^d$-graded algebra  over  an  artinian local ring $A$  and \; let  $M$ be a finitely generated $\mathbb{N}^d$-graded $S$-module such that  $S_{(1,1,\ldots,1)}$ is not contained in $ \sqrt{\mathrm{Ann} M}$.
Denote by  $\Lambda$ the set of all homogeneous prime ideals $P$ of  $S$ such that $P \in \mathrm{Supp}_{++}M$ and $\dim \mathrm{Proj}(S/P) = \dim \mathrm{Supp}_{++}M .$   Then  $$e(M;k_1,\ldots,k_d)= \sum_{P \in \Lambda}\ell(M_P)e(S/P;k_1,\ldots,k_d).$$}

We would like to emphasize that
although Theorem 3.1 is a general result for mixed multiplicities of  multi-graded modules which is a general object for mixed multiplicities of ideals, up to now we can not prove  the following  theorem  by using Theorem 3.1.

\vskip 0.2cm

\noindent{\bf Theorem 3.2.} {\it Let  $(R, \frak n)$  be  a  noetherian   local ring  with maximal ideal $\frak{n},$  infinite residue field $k = R/\frak{n},$ ideals $I_1,\ldots,I_d,$ an  $\frak n$-primary $J.$ Let $N$ be a finitely generated  $R$-module. Assume  that $I=I_1\cdots I_d$  is not contained in $ \sqrt{\mathrm{Ann}{N}}.$  Set $\overline{N}=\dfrac{N}{0_N: I^\infty}.$ Denote by  $\Pi$ the set of all prime ideals $\frak p $ of  $R$ such that $\frak p \in \mathrm{Min}(R/\mathrm{Ann}\overline{N})$ and $\dim R/\frak p = \dim \overline{N}.$
Then we have
   $$e(J^{[k_0+1]},I_1^{[k_1]},\ldots,I_d^{[k_d]}; N)= \sum_{\frak p \in \Pi}\ell({N}_{\frak p})e(J^{[k_0+1]},I_1^{[k_1]},\ldots,I_d^{[k_d]}; R/\frak p).$$}
\;\;\;It is natural to suppose that the proof of Theorem 3.2 will have to use the additive property on exact sequences of mixed multiplicities. But in fact, this approach seems to become a obstruction in proving  Theorem 3.2. This is a motivation to help us giving another approach for the proof of Theorem 3.2 as in this paper (see the proof of Theorem 3.2, Section 3).
On the contrary,  even  from  Theorem 3.2 we show that mixed multiplicities of arbitrary ideals are additive on exact sequences (see Corollary 3.9, Section 3) which covers \cite[Lemma 17.4.4]{HS}.

      Our approach  is based on
multiplicity formulas of multi-graded Rees modules with respect to powers of ideals (see Proposition  2.4 and Corollary 2.5, Section 2) via  linking minimal homogeneous prime ideals of maximal coheight of the Rees module $\mathfrak R(I_1,\ldots,I_d; N)= \bigoplus_{n_1,\ldots,n_d\ge 0}I_1^{n_1}\cdots I_d^{n_d}N$ and minimal prime ideals of maximal coheight of $N$ (see Lemma 3.4, Section 3).

Set
$$S_{\widehat{i}} \bigoplus_{n_1,\ldots,n_{i-1}, n_{i+1},\ldots, n_d \ge 0\;;\;n_i=0}S_{(n_1,\ldots,n_d)}\; \text {and}\; M_{\widehat{i}} = S_{\widehat{i}}M_{(0,\ldots,0)}.$$

Next, we establish the recursion formulas for the sum of all the  mixed multiplicities of the $\mathbb{N}^{d}$-graded  module $M:$ $\widetilde{e}(M)= \sum_{k_1\:+\:\cdots\:+\:k_d\;=\;\ell-1}e(M;k_1,\ldots,k_d)$ which  express $\widetilde{e}(M)$ as a sum  $$\widetilde{e}(M) = \widetilde{e}(M/xM) + \widetilde{e}(W),$$ where $$\dim \text{Supp}_{++}(M/xM) = \dim \text{Supp}_{++}M-1$$
 and $W$ is  an $\mathbb{N}^{d-1}$-graded  module. This result can be stated as follows.

\vskip 0.2cm

 \noindent
{\bf Theorem  5.2.}\; {\it Let $S$  be  a finitely generated standard\; $\mathbb{N}^d$-graded  algebra  over an artinian local ring $A$  and  let $M$ be a finitely generated $\mathbb{N}^d$-graded $S$-module  such that  $M=SM_{(0,\ldots,0)}.$  Set $\dim_{S^\triangle} M^\triangle = \ell.$  Assume that   $e(M;k_1,\ldots,k_d)\ne 0 $ and $k_i > 0.$  Let $x\in S_i$ be an $S_{++}$-filter-regular element with respect to  $M.$ Set  $\mathrm{\bf h}= h_1,\ldots,h_d$ and $\mid\mathrm{\bf h}\mid = h_1+\cdots+h_d.$ Then the following statements hold.

\begin{itemize}
\item[$\mathrm{(i)}$] $\widetilde{e}\Big(\dfrac{M}{xM}\Big) =  \sum_{\mid\mathrm{\bf h}\mid\:=\:\ell-1;\; h_i >0}e(M;\mathrm{\bf h}).$
\item[$\mathrm{(ii)}$] $\sum_{\mid \mathrm{\bf h}\mid\:=\:\ell-1;\; h_i =0}e(M;\mathrm{\bf h})\ne 0$ if and only if
$\dim_{{S_{\widehat{i}}}^\triangle}[S_i^vM_{\widehat{i}}]^\triangle =\ell$ for some $v \gg 0.$ In this case,
 $\widetilde{e}\big({S_i^vM_{\widehat{i}}}\big)=\sum_{
\mid\mathrm{\bf h}\mid\:=\:\ell-1;\; h_i =0}e(M;\mathrm{\bf h})$ for all $v \gg 0.$
\item[$\mathrm{(iii)}$] If\; $\dim_{{S_{\widehat{i}}}^\triangle}[S_i^vM_{\widehat{i}}]^\triangle =\ell$ for some $v \gg 0$ then
$\widetilde{e}(M)= \widetilde{e}\Big(\dfrac{M}{xM}\Big)+ \widetilde{e}\big({S_i^vM_{\widehat{i}}}\big)$ $ \text{ for all}\;\;  v \gg 0.$
\item[$\mathrm{(iv)}$] If\; $\dim_{{S_{\widehat{i}}}^\triangle}[S_i^vM_{\widehat{i}}]^\triangle <\ell$ for some $v \gg 0$ then $\widetilde{e}(M)= \widetilde{e}\Big(\dfrac{M}{xM}\Big)$
   for all $v \gg 0.$
 \end{itemize}}

As  consequences  of  Theorem 5.2,  we get the recursion formulas for the  multiplicity of  multi-graded Rees modules (see Theorem 5.5; Corollary 5.6;  Corollary 5.7 and Corollary 5.8,  Section 5).

\vskip 0.2cm
The main results of this paper yield many interesting consequences such as
the additivity and reduction formulas for  mixed multiplicities of ideals of positive height that covers \cite[Theorem 17.4.8]{HS} for the case  of $\frak n$-primary ideals; the additive property on exact sequences for mixed multiplicities of ideals and the multiplicity of  multi-graded Rees modules; the recursion formulas for the multiplicity of  multi-graded Rees modules; and the multiplicity formulas of Rees modules.

\vskip 0.2cm

This paper is divided into  five  sections. Section 2 is devoted to the discussion of mixed multiplicities of multi-graded Rees modules and the multiplicity of Rees modules with respect to powers of ideals  (Proposition  2.4 and Corollary 2.5) that will be used as a tool in the proofs of the paper. Section 3 gives the additivity and reduction formulas for  mixed multiplicities of multi-graded modules and mixed multiplicities of arbitrary ideals.  Section 4 investigates the
relationship between  filter-regular sequences of  multi-graded $F_J(J,I_1,\ldots,I_d;R)$-module $ F_J(J,I_1,\ldots,I_d;N)$ and  weak-(FC)-sequences of ideals that  will be used in the proofs of  Section 5 (Proposition  4.5).
Section 5 introduces the recursion formulas for  the sum of  all the mixed multiplicities of  multi-graded  modules. And as an application, we obtain
the recursion formulas for  the  multiplicity of  multi-graded Rees modules.

\vskip 0.2cm

\centerline{\Large\bf2.  Multiplicity  of multi-graded Rees  modules}$\\$
\noindent   This section studies  mixed multiplicities and the multiplicity of multi-graded modules. We  will give   multiplicity and mixed multiplicity  formulas of  Rees modules with respect to powers of ideals that will be used as a tool in the proofs of the paper.

 Set $\dim M^\triangle = \ell.$   By \cite[Theorem 4.1]{HHRT},
    $\ell_A[M_{(n_1,\ldots,n_d)}]$ is a polynomial of degree $\dim \text{Supp}_{++}M$ for all large $n_1,\ldots,n_d.$  Remember  that   $\dim \text{Supp}_{++}M=\ell-1$ by \cite[Remark 3.1]{VM}.
  The terms of total degree $\ell-1$ in this polynomial have the form
$$B_M(n_1,n_2,\ldots,n_d)= \sum_{k_1\:+\:\cdots\:+\:k_d\;=\;\ell-1}e(M;k_1,\ldots,k_d)\dfrac{n_1^{k_1}\cdots n_d^{k_d}}{k_1!\cdots k_d!}.$$

\noindent
Then $e(M;k_1,\ldots,k_d)$ are   non-negative integers not all zero, called the {\it  mixed multiplicity of type $(k_1,\ldots,k_d)$ of $M$} \cite{HHRT}.
And from now on $B_M(n_1,n_2,\ldots,n_d)$ is called the {\it   Bahattacharya homogeneous polynomial} of $M$ \cite{Bh}.

Set  $\mathrm{\bf k}= k_1,\ldots,k_d$ and $\mid\mathrm{\bf k}\mid = k_1+\cdots+k_d.$  Denote by $\widetilde{e}(M)$ the sum  of all the mixed multiplicities of $M,$ i.e.,
$\widetilde{e}(M):=\sum_{\mid \mathrm{\bf k}\mid=\:\ell-1}e(M;\mathrm{\bf k}).$ It is well known that in generally, the multiplicity $e(M)$ of $M$ and $\widetilde{e}(M)$ are different invariants of $M$.

 Let  $(R, \frak n)$  be  a  noetherian   local ring  with maximal ideal $\frak{n},$  infinite residue field $k = R/\frak{n}$ and  let $N$ be a finitely generated  $R$-module. Let $I_1,\ldots,I_d$ be ideals of $R$
  such that $I_1\cdots I_d$  is not contained in $ \sqrt{\mathrm{Ann}{N}}.$

 \noindent
 Put $\mathrm{\bf I}= I_1,\ldots,I_d;$ $\mathrm{\bf n}= n_1,\ldots,n_d;$
 $\mathbb{I}^{\mathrm{\bf n}}= I_1^{n_1},\ldots,I_d^{n_d};$
 $\mathrm{\bf I}^{[\mathrm{\bf k}]}= I_1^{[k_1]},\ldots,I_d^{[k_d]}.$

 \noindent
 Denote by $$\mathfrak R(\mathrm{\bf I}; R) = \mathfrak R(I_1,\ldots,I_d;R)= \bigoplus_{n_1,\ldots,n_d\ge 0}I_1^{n_1}\cdots I_d^{n_d}$$ the Rees algebra of ideals $I_1,\ldots,I_d$ and by $$\;\;\;\mathfrak R(\mathrm{\bf I}; N) = \mathfrak R(I_1,\ldots,I_d;N)= \bigoplus_{n_1,\ldots,n_d\ge 0}I_1^{n_1}\cdots I_d^{n_d}N$$ the Rees module of ideals $I_1,\ldots,I_d$ with respect to $N.$
Let $J$ be an  $\frak n$-primary ideal. Set
$$F_J(J,\mathrm{\bf I}; R)= F_J(J,I_1,\ldots,I_d; R) =\bigoplus_{n_0, n_1,\ldots,n_d\ge 0}\dfrac{J^{n_0}I_1^{n_1}\cdots I_d^{n_d}}{J^{n_0+1}I_1^{n_1}\cdots I_d^{n_d}}$$
and $$\;\;\;\;\;F_J(J,\mathrm{\bf I}; N)= F_J(J,I_1,\ldots,I_d;N) =\bigoplus_{n_0, n_1,\ldots,n_d\ge 0}\dfrac{J^{n_0}I_1^{n_1}\cdots I_d^{n_d}{N}}{J^{n_0+1}I_1^{n_1}\cdots I_d^{n_d}{N}}.$$
Then $F_J(J,\mathrm{\bf I};R)$ is a finitely generated standard multi-graded algebra over an  artinian local ring $R/J$ and $F_J(J,\mathrm{\bf I}; N)$ is a finitely generated multi-graded $F_J(J,\mathrm{\bf I}; R)$-module.
Set $ I = I_1\cdots I_d.$
Denote by $B_N\big(J, \mathrm{\bf I};n_0, \mathrm{\bf n}\big)= B_N\big(J, \mathrm{\bf I};n_0, n_1,\ldots, n_d\big)$ the  Bahattacharya homogeneous polynomial of $ F_J(J,\mathrm{\bf I}; N).$ Then remember that
$$\deg B_N\big(J, \mathrm{\bf I};n_0, \mathrm{\bf n}\big) =\dim \dfrac{N}{0_N: I^\infty}-1$$ by \cite[Proposition 3.1]{Vi1} (see \cite[Proposition 3.1]{MV}).
And by \cite[Remark 3.1]{VM},
$\deg B_N\big(J, \mathrm{\bf I};n_0,\mathrm{\bf n}\big)=  \dim F_J(J,\mathrm{\bf I}; N)^\triangle-1.$
Hence
$\dim F_J(J,\mathrm{\bf I}; N)^\triangle = \dim \dfrac{N}{0_N: I^\infty}.$
In the case that $\mathrm{ht} \dfrac{I+\mathrm{Ann}N}{\mathrm{Ann} N} > 0,$  $\dim \dfrac{N}{0_N: I^\infty} = \dim N.$
The above facts yield:

\vskip 0.2cm
\noindent {\bf Note 2.1.} $\dim F_J(J,\mathrm{\bf I}; N)^\triangle = \dim \dfrac{N}{0_N: I^\infty},$ \;and if $\;\mathrm{ht} \dfrac{I+\mathrm{Ann}N}{\mathrm{Ann} N} > 0$ then
$$\dim F_J(J,\mathrm{\bf I}; N)^\triangle = \dim {N}.$$
Set $\dim \dfrac{N}{0_N:{ I}^\infty} = q$ and
$$ e\big(F_J(J,\mathrm{\bf I}; N);k_0, k_1,\ldots,k_d\big)  = e\big(J^{[k_0+1]},I_1^{[k_1]},\ldots,I_d^{[k_d]};N\big):=  e\big(J^{[k_0+1]},\mathrm{\bf I}^{[\mathrm{\bf k}]};N\big)$$ $(k_0+k_1+\cdots+k_d = k_0 + \mid\mathrm{\bf k}\mid = q-1).$ Then $e\big(J^{[k_0+1]},I_1^{[k_1]},\ldots,I_d^{[k_d]};N\big) $ is called the   {\it mixed multiplicity of $N$ with respect to ideals $J,I_1,\ldots,I_d$
of type $(k_0,k_1,\ldots,k_d)$} (see \cite{MV, Ve}).

\vskip 0.2cm
\noindent {\bf Note 2.2.} Recall that \; by \cite[Proposition 3.1]{MV} which is a generalized result \;of \; \cite[Proposition 3.1]{Vi1},  we have
$e\big(J^{[k_0+1]},\mathrm{\bf I}^{[\mathrm{\bf k}]}; N\big)= e\Big(J^{[k_0+1]},\mathrm{\bf I}^{[\mathrm{\bf k}]};\dfrac{N}{0_N: I^\infty}\Big),$ and hence
$$\widetilde{e}\big(F_J(J,\mathrm{\bf I}; {N})\big)= \widetilde{e}\Big(F_J(J,\mathrm{\bf I}; \dfrac{N}{0_N: I^\infty})\Big).$$

\vskip 0.2cm
\noindent {\bf Note 2.3.} By \cite[Corollary 4.6]{HHRT}, it follows that
$$e\big(\big(J,\mathfrak R(\mathrm{\bf I}; R)_+\big); \mathfrak R(\mathrm{\bf I}; N)\big)= e\big(F_J(J,\mathrm{\bf I}; N)\big).$$

\noindent
Now, assume that $\mathrm{ht} \dfrac{I+\mathrm{Ann}N}{\mathrm{Ann} N} > 0.$ Then $\dim \dfrac{N}{0_N: I^\infty} = \dim N.$ In this case,

$$B_N\big(J, \mathrm{\bf I};n_0,\mathrm{\bf n}\big)
= \sum_{k_0\:+\mid\mathrm{\bf k}\mid\;=\;q-1}e\big(J^{[k_0+1]},\mathrm{\bf I}^{[\mathrm{\bf k}]};N\big)
\dfrac{n_0^{k_0}n_1^{k_1}\cdots n_d^{k_d}}{k_0!k_1!\cdots k_d!}\eqno(1)$$
and $$e\big(\big(J,\mathfrak R(\mathrm{\bf I}; R)_+\big); \mathfrak R(\mathrm{\bf I}; N)\big)= \sum_{k_0\:+\:\mid\mathrm{\bf k}\mid
=\;q-1}e\big(J^{[k_0+1]},\mathrm{\bf I}^{[\mathrm{\bf k}]};N\big)\eqno(2)$$
by \cite[Theorem 4.4]{HHRT} which is a generalized version of \cite[Theorem 1.4]{Ve}.
Next, let $u_1,\ldots,u_d$ be positive integers. Set $\mathrm{\bf u}^\mathrm{\bf k}= u_1^{k_1}\cdots u_d^{k_d}.$  From (1) we have
\begin{align*}
&B_N\big(J, \mathrm{\bf I}^\mathrm{\bf u}, n_0 , \mathrm{\bf n}\big)
= \sum_{k_0\:+\mid\mathrm{\bf k}\mid\;=\;q-1}e\big({J}^{[k_0+1]},{\mathrm{\bf I}^{\mathrm{\bf u}}}^{[\mathrm{\bf k}]};N\big)
\dfrac{n_0^{k_0}n_1^{k_1}\cdots n_d^{k_d}}{k_0!k_1!\cdots k_d!}\;\;\; \mathrm{and}\\
&B_N\big(J, \mathrm{\bf I}^\mathrm{\bf u}, n_0 , \mathrm{\bf n}\big)
= \sum_{k_0\:+\mid\mathrm{\bf k}\mid\;=\;q-1}e\big(J^{[k_0+1]},\mathrm{\bf I}^{[\mathrm{\bf k}]};N\big)
\dfrac{n_0^{k_0}(u_1n_1)^{k_1}\cdots (u_dn_d)^{k_d}}{k_0!k_1!\cdots k_d!}. \end{align*}
Consequently,
$e\big({J}^{[k_0+1]},{\mathrm{\bf I}^{\mathrm{\bf u}}}^{[\mathrm{\bf k}]};N\big)=
e\big(J^{[k_0+1]},\mathrm{\bf I}^{[\mathrm{\bf k}]}; N\big)\mathrm{\bf u}^\mathrm{\bf k}.$  Hence by (2),
$$e\big(\big(J,\mathfrak R(\mathrm{\bf I}^\mathrm{\bf u}; R)_+\big); \mathfrak R(\mathrm{\bf I}^\mathrm{\bf u}; N)\big)= \sum_{k_0\:+\mid\mathrm{\bf k}\mid\;
=\;q-1}e\big(J^{[k_0+1]},\mathrm{\bf I}^{[\mathrm{\bf k}]}; N\big)\mathrm{\bf u}^\mathrm{\bf k}.$$
We obtain the following result.
\vskip 0.2cm
\noindent{\bf Proposition  2.4.}\;{\it Assume that $\mathrm{ht} \dfrac{I+\mathrm{Ann}N}{\mathrm{Ann} N} > 0$ and  $u_1,\ldots,u_d$ are  positive integers.
   Then

 \begin{itemize}
 \item[$\mathrm{(i)}$]$e\big({J}^{[k_0+1]},{\mathrm{\bf I}^{\mathrm{\bf u}}}^{[\mathrm{\bf k}]};N\big)=
e\big(J^{[k_0+1]},\mathrm{\bf I}^{[\mathrm{\bf k}]}; N\big)\mathrm{\bf u}^\mathrm{\bf k}.$

\item[$\mathrm{(ii)}$]
$e\big(\big(J,\mathfrak R(\mathrm{\bf I}^\mathrm{\bf u}; R)_+\big); \mathfrak R(\mathrm{\bf I}^\mathrm{\bf u}; N)\big)= \sum_{k_0\:+\mid\mathrm{\bf k}\mid\;
=\;q-1}e\big(J^{[k_0+1]},\mathrm{\bf I}^{[\mathrm{\bf k}]}; N\big)\mathrm{\bf u}^\mathrm{\bf k}.$
\end{itemize}}
Set $\overline{N} = \dfrac{N}{0_N: I^\infty}.$  It can be verified that
$\mathrm{ht} \dfrac{I+\mathrm{Ann}\overline{N}}{\mathrm{Ann} \overline{N}} > 0.$ By Note 2.2,
$$e\big(J^{[k_0+1]},\mathrm{\bf I}^{[\mathrm{\bf k}]};N\big)= e\Big(J^{[k_0+1]},\mathrm{\bf I}^{[\mathrm{\bf k}]};\dfrac{N}{0_N: I^\infty}\Big).$$
Then as an immediate consequence of Proposition  2.4 we get the following.
\vskip 0.2cm
\noindent{\bf Corollary 2.5.}\;{\it Let $u_1,\ldots,u_d$ be positive integers.   Then
 \begin{itemize}
 \item[$\mathrm{(i)}$]$e\big({J}^{[k_0+1]},{\mathrm{\bf I}^{\mathrm{\bf u}}}^{[\mathrm{\bf k}]}; N\big)=
e\big(J^{[k_0+1]},\mathrm{\bf I}^{[\mathrm{\bf k}]}; N\big)\mathrm{\bf u}^\mathrm{\bf k}.$
\item[$\mathrm{(ii)}$]
$e\Big(\big(J,\mathfrak R(\mathrm{\bf I}^\mathrm{\bf u}; R)_+\big); \mathfrak R\big(\mathrm{\bf I}^\mathrm{\bf u}; \dfrac{N}{0_N: I^\infty}\big)\Big)= \sum_{k_0\:+\mid\mathrm{\bf k}\mid
=\;q-1}e\big(J^{[k_0+1]},\mathrm{\bf I}^{[\mathrm{\bf k}]};N\big)\mathrm{\bf u}^\mathrm{\bf k}.$
\end{itemize}}

\noindent
Set $\mathbb{S}= F_J(J,\mathrm{\bf I}; R)$
and $\mathbb{M}= F_J(J,\mathrm{\bf I};N).$
Recall that $$\widetilde{e}(\mathbb{M}) = \sum_{k_0\:+\:\mid\mathrm{\bf k}\mid
=\;q-1}e(J^{[k_0+1]},\mathrm{\bf I}^{[\mathrm{\bf k}]};N).$$
Hence combining this fact with  Note 2.3 and Corollary 2.5  yields:
\vskip 0.2cm
\noindent{\bf Corollary 2.6.}
$e\Big(F_J(J,\mathrm{\bf I}; \dfrac{N}{0_N: I^\infty})\Big)= \widetilde{e}(\mathbb{M})= e\Big(\big(J,\mathfrak R(\mathrm{\bf I}; R)_+\big); \mathfrak R(\mathrm{\bf I}; \dfrac{N}{0_N: I^\infty})\Big).$\\
\vskip 0.2cm
\noindent{\bf Remark 2.7.}
If $\mathrm{ht} \dfrac{I+\mathrm{Ann}N}{\mathrm{Ann} N} > 0,$ then
 $e\big(\mathbb{M}\big)$ is the sum  of all
the mixed multiplicities of $\mathbb{M}$ by \cite{ HHRT, Ve}. Hence
$e\big(\mathbb{M}\big)= \widetilde{e}\big(\mathbb{M}\big).$  Thus
$e\big(F_J(J,\mathrm{\bf I}; {N})\big)= e\Big(F_J(J,\mathrm{\bf I}; \dfrac{N}{0_N: I^\infty})\Big)$ and
$e\big(\big(J,\mathfrak R(\mathrm{\bf I}; R)_+\big); \mathfrak R(\mathrm{\bf I}; {N})\big)= e\Big(\big(J,\mathfrak R(\mathrm{\bf I}; R)_+\big); \mathfrak R(\mathrm{\bf I}; \dfrac{N}{0_N: I^\infty})\Big)$ by
Corollary 2.6.

\vspace*{12pt}

\centerline{\Large\bf3. Additivity and reduction formulas for  mixed multiplicities  }

\vspace*{12pt}

\noindent
In this section, we prove additivity and reduction formulas for  mixed multiplicities.
And as an application of these formulas, we  show that mixed multiplicities of arbitrary ideals are additive  on exact sequences.

First, we have the following result for $\mathbb{N}^d$-graded $S$-modules.
\vskip 0.2cm
\noindent
{\bf Theorem 3.1.}\;{\it Let $S$  be  a finitely generated standard \;$\mathbb{N}^d$-graded algebra  over  an  artinian local ring $A$  and \; let  $M$ be a finitely generated $\mathbb{N}^d$-graded $S$-module such that  $S_{(1,1,\ldots,1)}$ is not contained in $ \sqrt{\mathrm{Ann} M}$.
Denote by  $\Lambda$ the set of all homogeneous prime ideals $P$ of  $S$ such that $P \in \mathrm{Supp}_{++}M$ and $\dim \mathrm{Proj}(S/P) = \dim \mathrm{Supp}_{++}M .$   Then  $$e(M;\mathrm{\bf k})= \sum_{P \in \Lambda}\ell(M_P)e(S/P;\mathrm{\bf k}).$$}
\begin{proof}\;\;  Denote by $B_M(\mathrm{\bf n})$ the  Bahattacharya homogeneous polynomial of $M.$ Remember \; that \;since \; $S_{(1,1,\ldots,1)}\nsubseteq \; \sqrt{\mathrm{Ann} M},$ \;
$\deg B_M(\mathrm{\bf n}) = \dim\text{Supp}_{++}M $ \; by \cite[Theorem 4.1]{HHRT} (see \cite[Remark 3.1]{VM}).  Let $$0 = M_0 \subseteq M_1 \subseteq M_2 \subseteq\cdots \subseteq M_u=M$$ be a prime  filtration of M, i.e., $M_{i+1}/M_i \cong S/P_i$ where  $P_i$ is a homogeneous prime ideal for all $0 \le i \le u-1.$  Since $S_{(1,1,\ldots,1)}\nsubseteq  \sqrt{\mathrm{Ann} M},\;\emptyset \ne \Lambda \subseteq \mathrm{Min}(S/\mathrm{Ann}M)$  by \cite[Lemma 1.1]{HHRT}.
Consequently,  $\Lambda \subseteq \{P_0,P_1,\ldots,P_{u-1}\}.$
 Note that $$\{P_0,P_1,\ldots,P_{u-1}\} \subseteq \text{Supp}M.$$ Hence if
$P_i \notin \mathrm{Supp}_{++}M$ then
$P_i \supseteq S_{++}.$ \; In this case, $\Big(\dfrac{S}{P_i}\Big)_{\mathrm{\bf n}}= 0$ for all $\mathrm{\bf n} \gg 0$ by \cite[Proposition 2.7]{VM}. Therefore  $B_{S/P_i} (\mathrm{\bf n})= 0.$
If $\dim \mathrm{Proj}(S/P_i) < \dim \mathrm{Supp}_{++}M ,$  we have
$\deg B_{S/P_i} (\mathrm{\bf n}) = \dim \mathrm{Proj}(S/P_i) < \dim \mathrm{Supp}_{++}M$ by
  \cite[Theorem 4.1]{HHRT}.  From the above facts, it follows that
 $$\deg B_{S/P_i} (\mathrm{\bf n}) < \dim \mathrm{Supp}_{++}M$$ for all $P_i \notin \Lambda.$
Hence $ B_M(\mathrm{\bf n})$ is a sum of all the  $B_{S/P}(\mathrm{\bf n})$ for $P \in \Lambda,$ counted as many times as $S/P$ appears as some $\dfrac{M_{i+1}}{M_i}.$ This number is exactly the length of $M_P$\; because   $\Lambda \subseteq \mathrm{Min}(S/\mathrm{Ann}M).$ Therefore
$B_M(\mathrm{\bf n}) =  \sum_{P \in \Lambda}\ell(M_P)B_{S/P}(\mathrm{\bf n}).$ Set $\dim \mathrm{Supp}_{++}M = s.$
Remember that $\mathrm{\bf n}^\mathrm{\bf k}:= n_1^{k_1}\cdots n_d^{k_d}.$ Now since
$$\;\;B_{S/P}(\mathrm{\bf n})=\sum_{\mid \mathrm{\bf k}\mid =\;s}
e(S/P;\mathrm{\bf k})\dfrac{\mathrm{\bf n}^\mathrm{\bf k}}{k_1!\cdots k_d!} \;\;\mathrm {for \;\; any} \; P \in \Lambda,$$
$$B_M(\mathrm{\bf n})= \sum_{\mid \mathrm{\bf k}\mid\;=\;s}\Big[\sum_{P \in \Lambda}\ell(M_P)e(S/P;\mathrm{\bf k})\Big]\dfrac{\mathrm{\bf n}^\mathrm{\bf k}}{k_1!\cdots k_d!}.$$ Hence
  $$\begin{array}{l}\sum_{\mid\mathrm{\bf k}\mid\;
=\;s}e(M;\mathrm{\bf k})\dfrac{\mathrm{\bf n}^\mathrm{\bf k}}{k_1!\cdots k_d!}\\
= \sum_{\mid \mathrm{\bf k}\mid\;=\;s}\Big[\sum_{P \in \Lambda}\ell(M_P)e(S/P;\mathrm{\bf k})\Big]\dfrac{\mathrm{\bf n}^\mathrm{\bf k}}{k_1!\cdots k_d!}.\end{array}$$
Thus, $$e(M;\mathrm{\bf k})= \sum_{P \in \Lambda}\ell(M_P)e(S/P;\mathrm{\bf k}).\;\blacksquare$$
\end {proof}

Although  Theorem 3.1 is a general result for mixed multiplicities of  multi-graded modules that is a general object for mixed multiplicities of ideals,
up to now we can not give a proof for the case of mixed multiplicities of ideals in the following result  by  using this theorem.
\vskip 0.2cm
\noindent
{\bf Theorem 3.2.} {\it Let  $(R, \frak n)$  be  a  noetherian   local ring  with maximal ideal $\frak{n},$  infinite residue field $k = R/\frak{n},$ ideals $I_1,\ldots,I_d,$ an  $\frak n$-primary $J.$ Let $N$ be a finitely generated  $R$-module. Assume  that $I=I_1\cdots I_d$  is not contained in $ \sqrt{\mathrm{Ann}{N}}.$  Set $\overline{N}=\dfrac{N}{0_N: I^\infty}.$ Denote by  $\Pi$ the set of all prime ideals $\frak p $ of  $R$ such that $\frak p \in \mathrm{Min}(R/\mathrm{Ann}\overline{N})$ and $\dim R/\frak p = \dim \overline{N}.$
Then we have
   $$e(J^{[k_0+1]},\mathrm{\bf I}^{[\mathrm{\bf k}]};N)= \sum_{\frak p \in \Pi}\ell({N}_{\frak p})e(J^{[k_0+1]},\mathrm{\bf I}^{[\mathrm{\bf k}]};R/\frak p).$$}
\vskip 0.2cm
\noindent
{\bf Remark 3.3.} Recall that  $\Pi$ is the set of prime ideals $\frak p$ such that $\frak p \in \mathrm{Min}(R/\mathrm{Ann}\overline{N})$ and $\dim R/\frak p = \dim \overline{N}.$ It is easy seen that
$$\Pi = \Big\{\frak p \in \mathrm{Ass}\Big(\frac{R}{\mathrm{Ann}\overline{N}}\Big) \mid  \; \dim R/\frak p = \dim \overline{N} \Big\}. $$
Since
$\mathrm{Ann}\overline{N} = \mathrm{Ann}N:I^\infty,$ $\frac{R}
{\mathrm{Ann}\overline{N}} = \frac{R}{\mathrm{Ann}N:I^\infty}.$
Consequently
\begin{align*}
\Pi &= \Big\{\frak p \in \mathrm{Ass}\Big(\frac{R}{\mathrm{Ann}N:I^\infty}\Big) \mid  \; \dim R/\frak p = \dim \overline{N} \Big\}\\
&= \Big\{\frak p \in \mathrm{Ass}\Big(\frac{R}{\mathrm{Ann}N}\Big) \mid  \; \frak p \nsupseteq I \; \mathrm{and}\; \dim R/\frak p = \dim \overline{N} \Big\}.  \end{align*}
If $\frak p \in \Pi$, $\overline{N}_{\frak p} = N_{\frak p}$ because  $I \nsubseteq \frak p.$
 Since $\ell({N}_{\frak p})= \ell(\overline{N}_{\frak p}) < +\infty,$ $ \frak p \in \mathrm{Min}(\frac{R}{\mathrm{Ann}N}).$
Hence
$\Pi = \Big\{\frak p \in \mathrm{Min}\Big(\frac{R}{\mathrm{Ann}N}\Big) \mid  \; \frak p \nsupseteq I \; \mathrm{and}\; \dim R/\frak p = \dim \overline{N} \Big\}.$
In the case that $\mathrm{ht}\dfrac{I + \mathrm{Ann}N}{\mathrm{Ann}N}> 0,$  $\dim \overline{N} = \dim N$ and $\frak p \nsupseteq I$ for any
$\frak p \in \mathrm{Min}(\frac{R}{\mathrm{Ann}N}).$    Consequently
$\Pi = \Big\{\frak p \in \mathrm{Min}\Big(\frac{R}{\mathrm{Ann}N}\Big) \mid \;\; \dim R/\frak p = \dim N \Big\}. $

Our approach  is based on
multiplicity formulas of multi-graded Rees modules with respect to powers of ideals that gave in  Proposition  2.4 via  linking homogeneous prime ideals which are in $\mathrm{Min}(\mathfrak R(\mathbf{I}; R)/\mathrm{Ann}\;\mathfrak R(\mathbf{I}; N))$ of maximal coheight and prime ideals in $\Pi$  by  the following lemma.

\vskip 0.2cm

  \noindent
{\bf Lemma 3.4.} {\it Let $N$ be a finitely generated $R$-module and let  $I_1, \ldots, I_d$ be ideals of $R$ such that $\mathrm{ht}\dfrac{I + \mathrm{Ann}N}{\mathrm{Ann}N}> 0$ $(I = I_1\cdots I_d)$.
Denote by $\Lambda$  the set of homogeneous prime ideals $P$ of the Rees algebras $\mathfrak R(\mathbf{I}; R)$ such that
$P \in \mathrm{Min}(\mathfrak R(\mathbf{I}; R)/\mathrm{Ann}\mathfrak R(\mathbf{I}; N))$ and
$\dim \mathfrak R(\mathbf{I}; R)/P = \dim \mathfrak R(\mathbf{I}; N),$ and denote by  $\Pi$  the set of   prime ideals of $R$ such that $\frak p \in \mathrm{Min}(R/\mathrm{Ann}N)$ and  $\dim R/\frak p = \dim N$. Then
there is an one-to-one   correspondence
between the set of prime ideals \; $\Pi$  and the set of prime \;ideals  $\Lambda$ given by
$$\frak p \mapsto
P = \bigoplus_{n_1, \ldots n_d\geq 0}(\frak p \cap I_1^{n_1}\cdots I_d^{n_d}).$$}

 \begin{proof} \;\; First, remember that since $\mathrm{ht}\dfrac{I + \mathrm{Ann}N}{\mathrm{Ann}N}> 0,$ $\dim \mathfrak R(\mathbf{I}; N) = \dim N +d.$
 Note that $\Lambda \subseteq \mathrm{Ass}_{\mathfrak R(\mathbf{I}; R)}\mathfrak R(\mathbf{I}; N)$ and  $\Pi \subseteq  \mathrm{Ass}_RN$ and
$$\mathrm{Ann }\;\mathfrak R(\mathbf{I}; N) =   \bigoplus_{n_1, \ldots n_d\geq 0}(\mathrm{Ann}N \cap I_1^{n_1}\cdots I_d^{n_d}).$$
    Now, if $\frak p$ is an ideal in $\Pi$, then it can be verified that
 $$P = \bigoplus_{n_1, \ldots n_d\geq 0}(\frak p \cap I_1^{n_1}\cdots I_d^{n_d})$$
is a homogeneous prime ideal of $\mathfrak R(\mathbf{I}; R)$ and $\mathrm{Ann}\;\mathfrak R(\mathbf{I}; N) \subseteq P.$

\noindent {\bf Note 3.5.} If $\frak q$ is a prime ideal of $R$ and
$I \nsubseteq \frak q$ then $\dfrac{I+ \frak q}{\frak q}\ne 0.$ Since $R/ \frak q$ is an integral domain and  $\dfrac{I+ \frak q}{\frak q}\ne 0,$   $\mathrm{ht}\dfrac{I+ \frak q}{\frak q} > 0.$ Therefore for any $\frak p \in \Pi,$ \;$\mathrm{ht}\dfrac{I+\frak p}{\frak p} > 0$  because  $I \nsubseteq \frak p$ by Remark 3.3.

It is easily seen that
\begin{align*}
\mathfrak R(\mathbf{I}; R)/P &= \bigoplus_{n_1, \ldots n_d\geq 0}\frac{I_1^{n_1}\cdots I_d^{n_d}}{\frak p \cap I_1^{n_1}\cdots I_d^{n_d}}\\
&\cong \bigoplus_{n_1, \ldots n_d\geq 0}\frac{I_1^{n_1}\cdots I_d^{n_d}+\frak p}{\frak p }
= \mathfrak R(\mathbf{I}; R/\frak p).
\end{align*}
 Since $\mathrm{ht}\dfrac{I + \frak p }{\frak p}> 0$ by Note 3.5,
$ \dim \mathfrak R(\mathbf{I}; R/\frak p) = \dim R/\frak p + d.$
Since $\frak p \in \Pi,$ $\dim R/\frak p =\dim N.$ Hence  $\dim \mathfrak R(\mathbf{I}; R)/P   = \dim \mathfrak R(\mathbf{I}; N).$
So  $P \in \Lambda.$

Next, suppose that $P$ is an ideal in $\Lambda$. Then $P$ is an associated prime ideal of $\mathfrak R(\mathbf{I}; N)$. Hence $P$ is homogeneous and there is a homogeneous element $x \in \mathfrak R(\mathbf{I}; N)$ such that $P = 0:x$. Set $\frak p = P\cap R$. Then
$\frak p = \{a \in R \mid ax = 0\}$
and $\mathrm{Ann}N \subseteq \frak p .$
Set $P = \bigoplus_{n_1, \ldots, n_d \geq 0}P_{(n_1, \ldots, n_d)}$, we have
$$P_{(n_1, \ldots, n_d)} = \{ a \in I_1^{n_1}\cdots I_d^{n_d} \mid ax = 0\}.$$
It implies that $P_{(n_1, \ldots, n_d)} = \frak p \cap I_1^{n_1}\cdots I_d^{n_d}.$ Therefore $P$ has the form
$$P = \bigoplus_{n_1, \ldots n_d\geq 0}(\frak p \cap I_1^{n_1}\cdots I_d^{n_d}).$$
Consequently $\mathfrak R(\mathbf{I}; R)/P \cong  \mathfrak R(\mathbf{I}; R/\frak p).$
Since $P \in \Lambda$ and $\mathrm{ht}\dfrac{I + \mathrm{Ann}N}{\mathrm{Ann}N}> 0,$
$$\dim \mathfrak R(\mathbf{I}; R)/P =  \dim \mathfrak R(\mathbf{I}; N) = \dim N +d.$$
Note that
$$\dim \mathfrak R(\mathbf{I}; R/\frak p) \leqslant \dim R/\frak p +d.$$
 Consequently $\dim R/\frak p \geqslant \dim N.$ Hence since $\mathrm{Ann}N \subseteq \frak p,$ $\dim R/\frak p =\dim N.$
Thus,
$\frak p \in \Pi.$
The above facts follow that there is  a bijection
between the set $\Pi$ and the set $\Lambda$ given by $\frak p \mapsto
P = \bigoplus_{n_1, \ldots n_d\geq 0}(\frak p \cap I_1^{n_1}\cdots I_d^{n_d}).$ $\blacksquare$
\end{proof}

\noindent
{\bf The proof of Theorem 3.2:} Let $u_1, \ldots, u_d$ be positive integers. Remember that  $\mathbf{I^u} = I_1^{u_1}, \ldots, I_d^{u_d}$.
 Set $\overline{N} =  \dfrac{N}{0_N : I^\infty}$ and $q = \dim \overline{N}$. Denote by $\Lambda_\mathbf{u}$ the set of  homogeneous prime ideals $P$ of the Rees algebra $\mathfrak R(\mathbf{I^u}; R)= \mathfrak R(I_1^{u_1}, \ldots, I_d^{u_d}; R)$ such that $P \in \mathrm{Min}(\mathfrak R(\mathbf{I^u}; R)/\mathrm{Ann}\mathfrak R(\mathbf{I^u}; \overline{N}))$
and $\dim \mathfrak R(\mathbf{I^u}; R)/P = \dim \mathfrak R(\mathbf{I^u}; \overline{N})$.
Recall that
$$\Pi = \Big\{\frak p \in \mathrm{Min}\Big(\frac{R}{\mathrm{Ann}\overline{N}}\Big) \mid  \;\;\; \dim R/\frak p = \dim \overline{N} \Big\}.$$
 By \cite[Theorem 11.2.4]{HS}, we have
$$\begin{array}{l}e\big(\big(J, \mathfrak R(\mathbf{I^u}; R)_+\big); \mathfrak R(\mathbf{I^u}; \overline{N})\big)\\
= \sum_{P\in \Lambda_\mathbf{u}}\ell(\mathfrak R(\mathbf{I^u}; \overline{N})_P)e\big(\big(J, \mathfrak R(\mathbf{I^u}; R)_+\big); \mathfrak R(\mathbf{I^u}; R)/P\big).\end{array} \eqno(3)$$
Remember that  $\mathrm{ht}\dfrac{I+ \mathrm{Ann}\overline{N}}{\mathrm{Ann}\overline{N}} > 0.$  In this case, if  $P \in \Lambda_\mathbf{u}$ and  $\frak p = P \cap R$, we have $$P = \bigoplus_{n_1, \ldots n_d\geq 0}(\frak p \cap {(I_1^{u_1})}^{n_1}\cdots {(I_d^{u_d})}^{n_d})$$ and $\frak p \in \Pi$ by Lemma 3.4. 
Next, we prove that $\ell(\mathfrak R(\mathbf{I^u};\overline{N})_P) = \ell(\overline{N}_{\frak p}).$ 
Indeed, since  
$\mathfrak R(\mathbf{I^u};R)_P/P\mathfrak R(\mathbf{I^u};R)_P \cong (\mathfrak R(\mathbf{I^u};R)/P)_P \cong \mathfrak R(\mathbf{I^u};R/{\frak p})_P,$ it follows that    
$\mathfrak R(\mathbf{I^u};R/{\frak p})_P$ is a simple $\mathfrak R(\mathbf{I^u};R)_P$-module. Now assume that $\ell_{R_{\frak p}}(\overline{N}_{\frak p}) = t.$ Then there exists a sequence of submodules of the $R$-module $\overline{N}:$
$$\overline{N} = N_0 \supset N_1 \supset \cdots \supset N_t= \{0\}$$ such that $(N_i/N_{i+1})_{\frak p} \cong R_{\frak p}/ 
 {\frak p} R_{\frak p}$ $(0 \leq i \leq t-1).$  
It can be verified that 
$$\mathfrak R(\mathbf{I^u};N_i/N_{i+1})_P \cong   \mathfrak R(\mathbf{I^u};(N_i/N_{i+1})_{\frak p})_P \cong \mathfrak R(\mathbf{I^u};R_{\frak p}/ 
 {\frak p} R_{\frak p})_P \cong \mathfrak R(\mathbf{I^u};R/{\frak p})_P.$$ 
So 
$\mathfrak R(\mathbf{I^u};N_i/N_{i+1})_P$ is a simple $\mathfrak R(\mathbf{I^u};R)_P$-module $(0 \leq i \leq t-1).$
By the above facts, we get a composition series of the $\mathfrak R(\mathbf{I^u};R)_P$-module $\mathfrak R(\mathbf{I^u};\overline{N})_P:$      
$$R(\mathbf{I^u};\overline{N})_P = R(\mathbf{I^u};N_0)_P \supset R(\mathbf{I^u};N_1)_P \supset \cdots \supset R(\mathbf{I^u};N_t)_P= \{0\}.$$  
Consequently $\ell(\mathfrak R(\mathbf{I^u};\overline{N})_P) = \ell(\overline{N}_{\frak p}).$
Hence since
$$\mathfrak R(\mathbf{I^u}; R)/P \cong \mathfrak R(\mathbf{I^u}; R/\frak p) \quad
$$
and by (3) we obtain

$$e\big((J, \mathfrak R(\mathbf{I^u}; R)_+); \mathfrak R(\mathbf{I^u}; \overline{N})\big)
= \sum_{\frak p\in \Pi}\ell(\overline{N}_\frak p)e\big((J, \mathfrak R(\mathbf{I^u}; R/\frak p)_+); \mathfrak R(\mathbf{I^u}; R/\frak p)\big). \eqno(4)$$
 Recall that $\frak p \in \Pi,$ $\mathrm{ht}\dfrac{I+\frak p}{\frak p} > 0$ by Note 3.5. Hence by Corollary 2.5(ii) and Proposition  2.4, we respectively get

$$\begin{array}{l}e\big((J, \mathfrak R(\mathbf{I^u}; R)_+); \mathfrak R(\mathbf{I^u}; \overline{N})\big)
= \sum_{k_0+ \mid\mathbf k\mid = q-1}e(J^{[k_0+1]}, \mathbf{I}^{[{\mathbf k}]}; N)\mathbf{u^k}
\end{array}\eqno(5)$$
and

$$\begin{array}{l}e\big((J, \mathfrak R(\mathbf{I^u}; R/\frak p)_+); \mathfrak R(\mathbf{I^u}; R/\frak p)\big)
= \sum_{k_0+ \mid \mathbf {k}\mid = q-1}e(J^{[k_0+1]}, \mathbf{I}^{[{\mathbf k}]}; R/\frak p)\mathbf{u^k}. \end{array}\eqno(6)$$
From (4), (5) and (6), it follows that
\begin{align*}
& \sum_{k_0+\mid\mathbf k\mid = q-1}e(J^{[k_0+1]}, \mathbf{I}^{[{\mathbf k}]}; N)\mathbf{u^k} \\
&= \sum_{\frak p \in \Pi}\ell(\overline{N}_\frak p)\Big( \sum_{k_0+\mid\mathbf k\mid = q-1}e(J^{[k_0+1]}, \mathbf{I}^{[{\mathbf k}]}; R/\frak p)\mathbf{u^k}\Big)\\
&= \sum_{k_0+\mid\mathbf k\mid = q-1}\Big(\sum_{\frak p\in \Pi}\ell(\overline{N}_\frak p) e(J^{[k_0+1]}, \mathbf{I}^{[{\mathbf k}]}; R/\frak p)\Big)\mathbf{u^k}.
\end{align*}
Therefore
$$e(J^{[k_0+1]}, \mathbf{I}^{[{\mathbf k}]}; N) = \sum_{\frak p\in \Pi}\ell(\overline{N}_\frak p) e(J^{[k_0+1]}, \mathbf{I}^{[{\mathbf k}]}; R/\frak p).$$ Recall  that $\ell(\overline{N}_\frak p) = \ell({N}_\frak p)
$ by Remark 3.3.   Thus
$$e(J^{[k_0+1]}, \mathbf{I}^{[{\mathbf k}]}; N) = \sum_{\frak p\in \Pi}\ell({N}_\frak p) e(J^{[k_0+1]}, \mathbf{I}^{[{\mathbf k}]}; R/\frak p).
\;\blacksquare$$

Note that if $\mathrm{ht}\dfrac{I + \mathrm{Ann}N}{\mathrm{Ann}N}> 0$ then  $\Pi = \big\{\frak p \in \mathrm{Min}\Big(\frac{R}{\mathrm{Ann}N}\Big) \mid \dim R/\frak p = \dim N \big\}$
by Remark 3.3. Hence by Theorem 3.2,  we obtain the following result.
\vskip 0.2cm
\noindent
{\bf Corollary 3.6.} {\it Let  $(R, \frak n)$  be  a  noetherian   local ring  with maximal ideal $\frak{n}$ and infinite residue field $k = R/\frak{n},$ ideals
$I_1,\ldots,I_d$ and  an  $\frak n$-primary ideal $J.$  Let $N$ be a finitely generated  $R$-module.   Set  $I=I_1\cdots I_d.$   Assume that $\mathrm{ht}\dfrac{I + \mathrm{Ann}N}{\mathrm{Ann}N}> 0.$  Denote by  $\Pi$ the set of all  prime ideals $\frak p $ of  $R$ such that $\frak p \in \mathrm{Min}(R/\mathrm{Ann}N)$ and $\dim R/\frak p = \dim N.$
 Then we have
   $$e(J^{[k_0+1]},\mathrm{\bf I}^{[\mathrm{\bf k}]};N)= \sum_{\frak p \in \Pi}\ell(N_{\frak p})e(J^{[k_0+1]}, \mathrm{\bf I}^{[\mathrm{\bf k}]};R/\frak p).$$}

Let $I_1,\ldots, I_d$ be $\frak n$-primary ideals of $R$. Set $\dim N = q.$
Denote by $P(n_1,\ldots,n_d)$ the Hilbert-Samuel polynomial of the Hilbert-Samuel function $\ell_A\Big(\frac{N}{I_1^{n_1}\cdots I_d^{n_d}N}\Big).$ For any $1\leq i\leq d,$ denote by $Q_i(n_1,\ldots,n_d)$ the Hilbert-Samuel polynomial of the Hilbert-Samuel function $\ell_A\Big(\frac{I_1^{n_1}\cdots I_i^{n_i}\cdots I_d^{n_d}N}{I_1^{n_1}\cdots I_i^{n_i+1}\cdots I_d^{n_d}N}\Big).$
Then we have $\deg P(n_1,\ldots,n_d) = q$ and
$$P(n_1,\ldots,n_i+1,\ldots,n_d) -P(n_1,\ldots,n_i,\ldots,n_d) = Q_i(n_1,\ldots, n_i, \ldots,n_d).$$
Write the terms of total degree $q$ in $P(\mathrm{\bf n})$ in the form
$\sum_{\mid\mathrm{\bf k}\mid = q} e(\mathrm{\bf I}^{[\mathrm{\bf k}]}; N)\frac{\mathrm{\bf n}^\mathrm{\bf k}}{k_1! \cdots k_d!}.$
 Since $k_1+\cdots+k_d= \mid\mathrm{\bf k}\mid = q >0,$
there exists $1 \le j \le d$ such that $k_j>0.$ It is easy to check that $\dfrac{e(\mathrm{\bf I}^{[\mathrm{\bf k}]}; N)}{k_1!\cdots (k_j-1)!\cdots k_d!}n_1^{k_1}\cdots n_j^{k_j-1}\cdots n_d^{k_d}$ is a term of total degree $q-1$ in $Q_j(\mathrm{\bf n}).$  So $e(\mathrm{\bf I}^{[\mathrm{\bf k}]}; N)$  as in \cite{HS} is exactly
the mixed multiplicity of $N$ with respect to $(I_1,\ldots,I_j,\ldots, I_d)$ of the type $(k_1,\ldots,k_j,\ldots,k_d)$ defined in Section 2 with $I_j$  playing the role of $J$. Therefore, for any non-negative integers $k_1,\ldots,k_d$ with $k_1+\cdots+k_d = \;\mid\mathrm{\bf k}\mid \;= q$, one also calls $e(\mathrm{\bf I}^{[\mathrm{\bf k}]}; N)$ the mixed multiplicity of $N$  with respect to $(I_1,\ldots, I_d)$ of the type $(k_1,\ldots,k_d).$

Then as a consequence of Corollary 3.6, we get the following result.

\vskip 0.2cm

     \noindent
{\bf Corollary 3.7} \cite[Theorem 17.4.8]{HS}. {\it Let  $(R, \frak n)$  be  a  noetherian   local ring  with maximal ideal $\frak{n}$ and infinite residue field $k = R/\frak{n},$ and $\frak n$-primary ideals  $I_1,\ldots,I_d.$    Let $N$ be a finitely generated  $R$-module of Krull dimension $\dim N  >0$.  Denote by  $\Pi$ the set of all prime ideals $\frak p $ of  $R$ such that $\frak p \in \mathrm{Min}(R/\mathrm{Ann}N)$ and $\dim R/\frak p = \dim N.$ Assume that $k_1,\ldots, k_d$ are non-negative integers with $k_1+\cdots+k_d = \dim N.$  Then we have
   $$e(\mathrm{\bf I}^{[\mathrm{\bf k}]};N)= \sum_{\frak p \in \prod}\ell(N_{\frak p})e(\mathrm{\bf I}^{[\mathrm{\bf k}]};R/\frak p).$$}
\begin{proof}\; Since $\dim N  >0$ and $I = I_1\cdots I_d$ is an $\frak n$-primary ideal, $\mathrm{ht}\dfrac{I + \mathrm{Ann}N}{\mathrm{Ann}N}> 0.$ Hence the proof is immediate from Corollary 3.6. $\blacksquare$

\end{proof}
\vskip 0.2cm
\noindent
{\bf Remark 3.8.}
Let $ W_1, W_2, W_3$ be finitely generated  $R$-modules and let $I_1,\ldots,I_d$ be ideals of $R$ such that $I=I_1\cdots I_d \nsubseteq \sqrt{\mathrm{Ann}{W_i}}$ \;for all $i = 1, 2, 3.$
Let $$ 0\longrightarrow W_1 \longrightarrow W_3 \longrightarrow W_2\longrightarrow 0 $$ be a short exact sequence of $R$-modules.
 For any $i = 1, 2, 3,$ set $\overline{W}_i= \dfrac{W_i}{0_{W_i}: I^\infty}$ and $p_i = \dim \overline{W}_i.$
Denote by $\Pi_i$ the set of  prime ideals such that  $\frak p \in \mathrm{Min}(R/\mathrm{Ann}\overline{W}_i)$ and  $\dim R/\frak p = p_i$.
Set $\Omega = \Pi_1 \cup\Pi_2\cup\Pi_3$. For any $\frak p \in \Omega,$ we have always
 short exact sequences
$$ 0\longrightarrow (W_1)_\frak p \longrightarrow (W_3)_\frak p \longrightarrow (W_2)_\frak p\longrightarrow 0. $$
If \; $p_j < p_i\;$\; and \;$p_k = \{p_1,\; p_2,\;p_3\}\setminus \{p_i,\;p_j\}$  then \;for any $\frak p \in \Pi_i,$ we get
$\;\dim \overline{W}_j <   \dim R/\frak p$, and hence
$\frak p \nsupseteq \mathrm{Ann}\overline{W}_j$. In this case, $(\overline{W}_j)_\frak p = 0.$  By  Remark 3.3,  $(W_j)_\frak p = (\overline{W}_j)_\frak p.$ Hence $(W_j)_\frak p = 0$. Thus
$0 \ne (W_i)_\frak p = (W_k)_\frak p.$ This argument proves that if $p_j < p_i$ then $p_i = p_k$ and $\Pi_i = \Pi_k,$ moreover, $p_3 = \max\{p_1, p_2\}$.

Using Theorem 3.2, now we prove that the mixed multiplicities of arbitrary ideals are additive  on short exact sequences by the following result.

\vskip 0.2cm
\noindent
{\bf Corollary 3.9.} {\it Keep the notations as in Remark $3.8$.
Let $J$ be an $\frak n$-primary ideal.
   Set $ \frak J = (J,\mathfrak R(\mathrm{\bf I}; R)_+).$  Assume that $\dim \overline{W}_3= k_0+k_1+\cdots+k_d+1.$  Then the following statements hold.
 \begin{itemize}
 \item[$\mathrm{(i)}$] If $\dim  \overline{W}_1 =\dim \overline{W}_2=\dim \overline{W}_3$
 then
\begin{align*}
(a )&: e(J^{[k_0+1]},\mathrm{\bf I}^{[\mathrm{\bf k}]}; W_3)=
e(J^{[k_0+1]},\mathrm{\bf I}^{[\mathrm{\bf k}]}; W_1)+
e(J^{[k_0+1]},\mathrm{\bf I}^{[\mathrm{\bf k}]}; W_2);\\
(b)&: e\big(\frak J; \mathfrak R(\mathrm{\bf I}; \overline{W}_3)\big)
= e\big(\frak J; \mathfrak R(\mathrm{\bf I}; \overline{W}_1)\big)+ e\big(\frak J; \mathfrak R(\mathrm{\bf I}; \overline{W}_2)\big).
\end{align*}

\item[$\mathrm{(ii)}$]  If $h\ne k = 1,2$ and $\dim \overline{W}_3 > \dim \overline{W}_h$ then
\begin{align*}
(a)&: e(J^{[k_0+1]},\mathrm{\bf I}^{[\mathrm{\bf k}]}; W_3)=
e(J^{[k_0+1]},\mathrm{\bf I}^{[\mathrm{\bf k}]}; W_k);\;\;\;\;\;\;\;\;\;\;\;\;\;\;\;\;\;\;\;\;\;\;\;\;\;\;\;\;\;\;\;\\
(b)&: e\big(\frak J; \mathfrak R(\mathrm{\bf I}; \overline{W}_3)\big)
= e\big(\frak J; \mathfrak R(\mathrm{\bf I}; \overline{W}_k)\big).
\end{align*}
\end{itemize}}

\begin{proof}\;\;
 The proof of (i): Since $p_1 = p_2 = p_3,$ by Theorem 3.2  we have
$$
e(J^{[k_0+1]}, \mathbf{I}^{[\mathbf{k}]}; W_i)
= \sum_{\frak p\in \Pi_i}\ell(W_i)_\frak p e(J^{[k_0+1]}, \mathbf{I}^{[\mathbf{k}]}; R/\frak p)
$$ for  $i = 1, 2, 3.$ Let $\frak p \in \Omega\setminus \Pi_i.$  Since $\dim R/\frak p = p_i,$  $\frak p \nsupseteq \mathrm{Ann}\overline{W}_i.$ Consequently, $(\overline{W}_i)_\frak p = 0.$  By Remark 3.3,
$(W_i)_\frak p = (\overline{W}_i)_\frak p.$ So $(W_i)_\frak p = 0.$
 From this it follows that
$$
e(J^{[k_0+1]}, \mathbf{I}^{[\mathbf{k}]}; W_i) = \sum_{\frak p\in \Pi_i}\ell(W_i)_\frak p e(J^{[k_0+1]}, \mathbf{I}^{[\mathbf{k}]}; R/\frak p)
= \sum_{\frak p\in \Omega}\ell(W_i)_\frak p e(J^{[k_0+1]}, \mathbf{I}^{[\mathbf{k}]}; R/\frak p)
$$ for all $i = 1, 2, 3.$
Therefore  by Theorem 3.2, we obtain
\begin{align*}
e(J^{[k_0+1]}, \mathbf{I}^{[\mathbf{k}]}; W_3) &= \sum_{\frak p\in \Omega}\ell(W_3)_\frak p e(J^{[k_0+1]}, \mathbf{I}^{[\mathbf{k}]}; R/\frak p)\\
&= \sum_{\frak p\in \Omega}(\ell(W_1)_\frak p + \ell(W_2)_\frak p) e(J^{[k_0+1]}, \mathbf{I}^{[\mathbf{k}]}; R/\frak p)\\
&= \sum_{\frak p\in \Omega}\ell(W_1)_\frak p e(J^{[k_0+1]}, \mathbf{I}^{[\mathbf{k}]}; R/\frak p)
+ \sum_{\frak p\in \Omega}\ell(W_2)_\frak p e(J^{[k_0+1]}, \mathbf{I}^{[\mathbf{k}]}; R/\frak p)\\
& = e(J^{[k_0+1]}, \mathbf{I}^{[\mathbf{k}]}; W_1) + e(J^{[k_0+1]}, \mathbf{I}^{[\mathbf{k}]}; W_2).
\end{align*}
Hence we get (a) of (i). By Corollary 2.5(ii) we have (b) of (i).
   The case that $p_3 > p_h:$ By Remark 3.8, $p_3 = p_k;$
 $\Pi_3 = \Pi_k$ and $(W_3)_\frak p = (W_k)_\frak p$ for all
$\frak p \in \Pi_3 = \Pi_k.$
 Consequently,  we obtain (ii) by Corollary 2.5(ii) and since
\begin{align*}
e(J^{[k_0+1]}, \mathbf{I}^{[\mathbf{k}]}; W_3) &= \sum_{\frak p\in \Pi_3}\ell(W_3)_\frak p e(J^{[k_0+1]}, \mathbf{I}^{[\mathbf{k}]}; R/\frak p)\\
&= \sum_{\frak p\in \Pi_k}\ell(W_k)_\frak p e(J^{[k_0+1]}, \mathbf{I}^{[\mathbf{k}]}; R/\frak p)\\
&= e(J^{[k_0+1]}, \mathbf{I}^{[\mathbf{k}]}; W_k).
 \; \blacksquare \end{align*}
\end{proof}

\noindent
{\bf Remark 3.10.}  Now, if we assign the mixed multiplicities of modules $W_i:$ $$e(J^{[k_0+1]},\mathrm{\bf I}^{[\mathrm{\bf k}]}; W_i) = 0$$ to  the case that $ k_0+\cdots+k_d > \dim \overline{W}_i-1,$ then from Corollary 3.9, we immediately get that: if  $k_0 + \mid \mathrm{\bf k}\mid = \dim\overline{W}_3-1$ then
$$e(J^{[k_0+1]},\mathrm{\bf I}^{[\mathrm{\bf k}]}; W_3)=
e(J^{[k_0+1]},\mathrm{\bf I}^{[\mathrm{\bf k}]}; W_1)+
e(J^{[k_0+1]},\mathrm{\bf I}^{[\mathrm{\bf k}]}; W_2).$$
 It is natural to suppose that the proof of Theorem 3.2 will be  based on Corollary 3.9. Hence
 one of obstructions in proving  Theorem 3.2 is Corollary 3.9. This is a motivation to help us giving the proof of Theorem 3.2 as in this paper.\\
 \newpage

  \vspace*{12pt}
\centerline{\Large\bf4. Filter-regular sequences of multi-graded modules  }
\vskip 0.2cm
\vspace*{12pt}
\noindent In this section, we
 explore  the
relationship between   filter-regular sequences of the multi-graded $F_J(J,\mathrm{\bf I}; R)$-module $ F_J(J,\mathrm{\bf I}; N)$   and  weak-(FC)-sequences of ideals that  will be used in the proofs of  Section 5.
\vskip 0.2cm
The concept of filter-regular sequences was introduced by Stuckrad and Vogel in \cite{SV}. The theory of filter-regular sequences became an important tool to study some classes of singular rings and has been continually developed (see e.g. \cite{BS, Hy,Tr1,  Tr2, VM}).
\vskip 0.2cm
\noindent{\bf Definition 4.1.} Let $S=\bigoplus_{n_1,\ldots,n_d\ge 0}S_{(n_1,\ldots,n_d)}$  be a finitely generated standard $\mathbb{N}^d$-graded algebra  over an artinian local ring $A$ and let $M=\bigoplus_{n_1,\ldots,n_d\ge 0}M_{(n_1,\ldots,n_d)}$ be  a finitely generated $\mathbb{N}^d$-graded $S$-module.
 Let $S_{(1,1,\ldots,1)}$ be not contained in $ \sqrt{\mathrm{Ann}M}$. Then
a homogeneous element $x\in S$ is called an {\it $S_{++}$-filter-regular element with respect to  $M$}  if $(0_M:x)_{(n_1,\ldots,n_d)}=0$ for  all large $n_1,\ldots,n_d.$
Let  $x_1,\ldots, x_t$ be homogeneous elements in $S$. We call that $x_1,\ldots, x_t$ is an {\it $S_{++}$-filter-regular sequence  with respect to  $M$} if
 $x_i$ is an $S_{++}$-filter-regular element with respect to $\dfrac{M}{(x_1,\ldots, x_{i-1})M}$ for all $i = 1,\ldots, t.$
\vskip 0.2cm
\noindent
{\bf Remark 4.2.} If  $S_{(1,1,\ldots,1)}\subseteq \sqrt{\mathrm{Ann} M}$ then $(0_M:x)_{(n_1,\ldots,n_d)} \subseteq M_{(n_1,\ldots,n_d)}= 0$ for  all large $n_1,\ldots,n_d.$  Hence any homogeneous element of $S$ always has the property of  an $S_{++}$-filter-regular element.
This fact only obstruct and do not carry usefulness. That is why in Definition 4.1, one has to exclude the case that $S_{(1,1,\ldots,1)}\subseteq \sqrt{\mathrm{Ann} M}$ in defining $S_{++}$-filter-regular elements.

\vskip 0.2cm
%\vspace*{12pt}
\noindent
{\bf Note 2.3.} If $S_{(1,1,\ldots,1)}\nsubseteq \sqrt{\mathrm{Ann} M},$ then by \cite{VM},
a homogeneous element $x\in S$ is an $S_{++}$-filter-regular element with respect to  $M$  if and only if
$x\notin P$ for any $P\in \mathrm{Ass}_SM$ and $P$ does not contain   $S_{++}.$ That means
$x\notin\bigcup_{S_{++}\nsubseteq P,\; P\in \mathrm{Ass}_SM}P.$
In this case, for any $1\le i \le d,$ there exists an $S_{++}$-filter-regular element $x \in S_i \setminus\frak m S_i.$

Remember that the positivity and the
relationship between  mixed multiplicities  and  Hilbert-Samuel multiplicities of ideals have attracted much attention (see e.g. \cite{KV, KR1, KR2, MV, Ro, Sw,Tr2,  Vi1, Vi2, Vi3,  VT}).
 In past years, using different sequences, one expressed
  mixed multiplicities into  Hilbert-Samuel multiplicity, for instance: Risler-Teissier in 1973 \cite{Te} by superficial sequences and Rees in 1984 \cite{Re} by joint reductions;  Viet in 2000 \cite{Vi1} by (FC)-sequences (see e.g. \cite{DV, MV,  VT}).

\vskip 0.2cm
\noindent {\bf Definition 4.4} \cite{Vi1}.  Let  $(R, \frak n)$ \; be  a  noetherian   local ring  with maximal ideal $\frak{n},$  infinite residue field $k = R/\frak{n}$ and let $N$ be a finitely generated  $R$-module. Let
$I_1,\ldots,I_d$ be ideals such that $I_1\cdots I_d$  is not contained in $ \sqrt{\mathrm{Ann}{N}}.$ Set $I=I_1\cdots I_d.$  An element $x \in R$ is called an {\it $(FC)$-element of $N$ with respect to $(I_1,\ldots, I_d)$}  if there exists $i \in \{ 1, \ldots, d\}$ such that $x \in I_i$ and
the following conditions are satisfied:

 \begin{itemize}
 \item[$\mathrm{(i)}$] $x$ is an $I$-filter-regular element with respect to $N,$ i.e.,\;$0_N:x \subseteq 0_N: I^{\infty}.$

\item[$\mathrm{(ii)}$] $x{N}\bigcap {I_1}^{n_1} \cdots I_i^{n_i+1}\cdots I_d^{n_d}{N}
= x{I_1}^{n_1}\cdots I_i^{n_i}\cdots I_d^{n_d}{N}$
for all $n_1,\ldots,n_d\gg0.$
\item [$\mathrm{(iii)}$] $\dim N/(xN:I^\infty)=\dim N/0_N:I^\infty-1.$
 \end{itemize}
       We call $x$ a {\it weak-$(FC)$-element of $N$ with respect to $(I_1,\ldots, I_d)$} if $x$ satisfies the conditions (i) and (ii).

Let $x_1, \ldots, x_t$ be a sequence in $R$. For any $0\le i < t,$ set\;      ${N}_i = \dfrac{N}{(x_1, \ldots, x_{i})N}$.  Then
$x_1, \ldots, x_t$ is called a {\it weak-$(FC)$-sequence of $N$ with respect to $(I_1,\ldots, I_d)$} if $x_{i + 1}$ is a weak-(FC)-element of ${N}_i$ with respect to $(I_1,\ldots, I_d)$  for all $i = 0, \ldots, t - 1$.

$x_1, \ldots, x_t$ is called an {\it $(FC)$-sequence of $N$ with respect to $(I_1,\ldots, I_d)$} if $x_{i + 1}$ is an (FC)-element of ${N}_i$ with respect to $(I_1,\ldots, I_d)$  for all $i = 0, \ldots, t - 1$.

 Recall that $$\begin{array}{l}\widetilde{e}(M)=\sum_{\mid \mathrm{\bf k}\mid=\:\ell-1}e(M;\mathrm{\bf k});\\   S_i=S_{(0,\ldots,{\underbrace{1}_i},\ldots,0)}
\mathrm {\;\;for \;\; all\;\;} i=1,\ldots,d;\\  \;\mathbb{S}= F_J(J,\mathrm{\bf I}; R) =\bigoplus_{n_0, n_1,\ldots,n_d\ge 0}\dfrac{J^{n_0}I_1^{n_1}\cdots I_d^{n_d}}{J^{n_0+1}I_1^{n_1}\cdots I_d^{n_d}};\\
\mathbb{M}= F_J(J,\mathrm{\bf I};N) =\bigoplus_{n_0, n_1,\ldots,n_d\ge 0}\dfrac{J^{n_0}I_1^{n_1}\cdots I_d^{n_d}{N}}{J^{n_0+1}I_1^{n_1}\cdots I_d^{n_d}{N}}.\end{array}$$
Set\;\;$S_{\widehat{i}} = \bigoplus_{n_1,\ldots,n_{i-1}, n_{i+1},\ldots, n_d \ge 0\;;\;n_i=0}S_{(n_1,\ldots,n_d)}\; \text {and}\; M_{\widehat{i}} = S_{\widehat{i}}M_{(0,\ldots,0)}.$

\vskip 0.2cm
\hspace*{-0.6cm}{\bf Proposition  4.5.}\;{\it Let $x \in I_i$ be a weak-$(FC)$-element  with respect to $(J, I_1,\ldots, I_d)$ of $N.$ Denote by
$\bar x$ and $\bar I_i$ the images  of $x$ and $I_i$ in $\mathbb{S}_i,$  respectively.  Then we have}:
\begin{itemize}
\item[$\mathrm{(i)}$]\;$\bar x$ is an {\it$\mathbb{S}_{++}$-filter-regular element with respect to $\mathbb{M}.$
\item[$\mathrm{(ii)}$]\;
$\dim (\mathbb{M}/\bar x\mathbb{M})^\triangle = \dim \dfrac{{N}}{x{N}: I^\infty}$ and $\widetilde{e}(\mathbb{M}/\bar x\mathbb{M}) =\widetilde{e}\big(F_J(J,\mathrm{\bf I};\dfrac{{N}}{x{N}})\big).$

\item[$\mathrm{(iii)}$]
$\begin{array}{l} \mathbb{S}/\bar I_i\mathbb{S}\cong F_J(J,I_1,\ldots,I_{i-1},I_{i+1},\ldots,I_d;R)\cong \mathbb{S}_{\widehat{i}}\;\text{and}\\
\mathbb{M}/\bar I_i\mathbb{M}\cong F_J(J,I_1,\ldots,I_{i-1},I_{i+1},\ldots,I_d;N)\cong\mathbb{M}_{\widehat{i}}.\end{array} $ }
\end{itemize}
\begin{proof}\; We have
$(0_{N}: I^\infty)\bigcap J^mI_1^{m_1}\cdots I_d^{m_d}{N}= 0$ for all  $m, m_1,\ldots, m_d \gg 0$ by Artin-Rees lemma. Since $x$ is an $ I$-filter-regular element with respect to $M,$ $$(0_N:x) \bigcap J^mI_1^{m_1}\cdots I_d^{m_d}{N} \subseteq (0_N: {I}^{\infty})\bigcap J^mI_1^{m_1}\cdots I_d^{m_d}{N}= 0$$
for all  $m, m_1,\ldots, m_d \gg 0.$
From this it follows that
\vspace*{12pt}
 \begin{align*}
 &\big(J^{n+1}I_1^{n_1}\cdots I_i^{n_i+1}\cdots I_d^{n_d}N:x\big)\bigcap J^{n}I_1^{n_1}\cdots I_i^{n_i}\cdots I_d^{n_d}N\\
&= \Big[\big(x{N }\bigcap J^{n+1}{I_1}^{n_1} \cdots I_i^{n_i+1}\cdots I_d^{n_d}{N }\big):x\Big]\bigcap J^{n}I_1^{n_1}\cdots I_i^{n_i}\cdots I_d^{n_d}N\\
 &= \Big[xJ^{n+1}{I_1}^{n_1} \cdots I_i^{n_i}\cdots I_d^{n_d}{N}:x\Big]\bigcap J^{n}I_1^{n_1}\cdots I_i^{n_i}\cdots I_d^{n_d}N\\
&= \Big[J^{n+1}{I_1}^{n_1} \cdots I_i^{n_i}\cdots I_d^{n_d}{N}+0_{N}:x\Big]\bigcap J^{n}I_1^{n_1}\cdots I_i^{n_i}\cdots I_d^{n_d}N\\
&= J^{n+1}{I_1}^{n_1} \cdots I_i^{n_i}\cdots I_d^{n_d}{N}+\big(0_{N}:x\big)\bigcap J^{n}I_1^{n_1}\cdots I_i^{n_i}\cdots I_d^{n_d}N\\
&= J^{n+1}{I_1}^{n_1} \cdots I_i^{n_i}\cdots I_d^{n_d}{N}
\end{align*}  for all $n, n_1,\ldots,n_d\gg0.$
\vspace*{10pt}
Hence $[0_{\mathbb{M}}: \bar x]_{(n, n_1,\ldots,n_d)} = 0$
for all $n, n_1,\ldots,n_d\gg0.$ Thus,  $\bar x$ is an $\mathbb{S}_{++}$-filter-regular element. We get (i). It can be verified that

$[\mathbb{M}/\bar x\mathbb{M}]_{(m, m_1,\ldots,m_d)}\cong\dfrac{J^mI_1^{m_1}\cdots I_d^{m_d}{N}}{J^{m+1}I_1^{m_1}\cdots I_d^{m_d}{N}+xJ^{m}I_1^{m_1}\cdots I_i^{m_i-1}\cdots I_d^{m_d}{N}} \;\;\; \text{and } $
 $$\begin{array}{l}\Big[F_J(J,\mathrm{\bf I};\dfrac{{N}}{x{N}})\Big]_{(m, m_1,\ldots,m_d)}=\bigg[\bigoplus_{n, n_1,\ldots,n_d\ge 0}\dfrac{J^nI_1^{n_1}\cdots I_d^{n_d}({N}/x{N})}{J^{n+1}I_1^{n_1}\cdots I_d^{n_d}({N}/x{N})}\bigg]_{(m, m_1,\ldots,m_d)}\\\cong \dfrac{J^mI_1^{m_1}\cdots I_d^{m_d}{N}+ x{N}}{J^{m+1}I_1^{m_1}\cdots I_d^{m_d}{N}+x{N}}
\vspace{6pt}\cong \dfrac{J^mI_1^{m_1}\cdots I_d^{m_d}{N}}{J^{m+1}I_1^{m_1}\cdots I_d^{m_d}{N} + x{N}\bigcap J^mI_1^{m_1}\cdots I_d^{m_d}{N}}.\end{array}$$
Since $x$ is a weak-(FC)-element,
$$x{N}\bigcap J^mI_1^{m_1}\cdots I_d^{m_d}{N}= xJ^{m}I_1^{m_1}\cdots I_i^{m_i-1}\cdots I_d^{m_d}{N}$$
for all $m, m_1,\ldots,m_d\gg0.$
Hence  $[\mathbb{M}/\bar{x}\mathbb{M}]_{(m, m_1,\ldots,m_d)} \cong \Big[F_J(J,\mathrm{\bf I};\dfrac{{N}}{x{N}})\Big]_{(m, m_1,\ldots,m_d)}$ for all $m, m_1,\ldots,m_d\gg0.$ From this it follows that $$\dim \big(\mathbb{M}/\bar x\mathbb{M}\big)^\triangle = \dim \Big[F_J(J,\mathrm{\bf I};\dfrac{{N}}{x{N}})\Big]^\triangle= \dim \dfrac{{N}}{x{N}: { I}^\infty}$$ by  Note 2.1  and $\widetilde{e}(\mathbb{M}/\bar x\mathbb{M}) =\widetilde{e}\Big(F_J(J,\mathrm{\bf I};\dfrac{{N}}{x{N}})\Big).$
We get (ii). Since $\bar I_i= \mathbb{S}_i,$ (iii) is obvious.
$\blacksquare$ \end{proof}
\vspace*{12pt}
\centerline{\Large\bf5. Recursion formulas for  multiplicities of
graded modules}
\vspace*{12pt}

\noindent   This section  gives the recursion formulas for the sum of all the  mixed multiplicities of multi-graded modules. And as an application, we get
the recursion formulas for  the  multiplicity of  multi-graded Rees modules.
Recall that  $\widetilde{e}(M)$ denotes  the sum  of all the mixed multiplicities of $M,$ i.e.,
$\widetilde{e}(M)=\sum_{\mid \mathrm{\bf k}\mid=\:\ell-1}e(M;\mathrm{\bf k});$
\begin{align*}
S_i&=S_{(0,\ldots,{\underbrace{1}_i},\ldots,0)}
\;\;\text {for all}\;\; i=1,\ldots,d; \\
S_{\widehat{i}} &= \bigoplus_{n_1,\ldots,n_{i-1}, n_{i+1},\ldots, n_d \ge 0\;;\;n_i=0}S_{(n_1,\ldots,n_d)}\; \text {and}\; M_{\widehat{i}} = S_{\widehat{i}}M_{(0,\ldots,0)}.
\end{align*}

   We have  the following comment.
\vskip 0.2cm
\noindent {\bf Remark 5.1.} For any $m \geqslant 0,$ $S_i^{m} M_{\widehat{i}}$ is a finitely generated $\mathbb{N}^{d-1}$-graded  $S_{\widehat{i}}\,$-module.
 Since  $0 : S_i^uM_{\widehat{i}} = 0 : S_i^{v} M_{\widehat{i}}$   for all $u, v \gg 0,$ there exists $h$ such that $\dim \mathrm{Supp}_{++}S_i^uM_{\widehat{i}}= \dim \mathrm{Supp}_{++}S_i^vM_{\widehat{i}}$ for all $u, v \ge h.$ Hence by \cite[Remark 3.1]{VM},
$\dim_{{S_{\widehat{i}}}^\triangle}[S_i^uM_{\widehat{i}}]^\triangle = \dim_{{S_{\widehat{i}}}^\triangle}[S_i^vM_{\widehat{i}}]^\triangle$ for all $u, v \ge h.$\\

The main result of this section is the following theorem.\\

 \noindent
{\bf Theorem  5.2.}\; {\it Let $S$  be  a finitely generated standard\; $\mathbb{N}^d$-graded  algebra  over an artinian local ring $A$  and  let $M$ be a finitely generated $\mathbb{N}^d$-graded $S$-module  such that  $M=SM_{(0,\ldots,0)}.$  Set $\dim_{S^\triangle} M^\triangle = \ell.$  Assume that   $e(M;k_1,\ldots,k_d)\ne 0 $ and $k_i > 0.$  Let $x\in S_i$ be an $S_{++}$-filter-regular element with respect to  $M.$  Then the following statements hold.

 \begin{itemize}
\item[$\mathrm{(i)}$]
 $\widetilde{e}\Big(\dfrac{M}{xM}\Big) = \sum_{\mid\mathrm{\bf h}\mid\:=\:\ell-1;\; h_i >0}e(M;\mathrm{\bf h}).$
\item[$\mathrm{(ii)}$] $\sum_{\mid \mathrm{\bf h}\mid\:=\:\ell-1;\; h_i =0}e(M;\mathrm{\bf h})\ne 0$ if and only if
$\dim_{{S_{\widehat{i}}}^\triangle}[S_i^vM_{\widehat{i}}]^\triangle =\ell$ for some $v \gg 0.$ In this case,
 $\widetilde{e}\big({S_i^vM_{\widehat{i}}}\big)=\sum_{
\mid\mathrm{\bf h}\mid\:=\:\ell-1;\; h_i =0}e(M;\mathrm{\bf h})$ for all $v \gg 0.$
\item[$\mathrm{(iii)}$] If\; $\dim_{{S_{\widehat{i}}}^\triangle}[S_i^vM_{\widehat{i}}]^\triangle =\ell$ for some $v \gg 0$ then
$\widetilde{e}(M)= \widetilde{e}\Big(\dfrac{M}{xM}\Big)+ \widetilde{e}\big({S_i^vM_{\widehat{i}}}\big)$ $ \text{ for all}\;\;  v \gg 0.$
\item[$\mathrm{(iv)}$] If\; $\dim_{{S_{\widehat{i}}}^\triangle}[S_i^vM_{\widehat{i}}]^\triangle <\ell$ for some $v \gg 0$ then $\widetilde{e}(M)= \widetilde{e}\Big(\dfrac{M}{xM}\Big)$
   for all $v \gg 0.$
 \end{itemize}}

\begin{proof}\;
Since $x\in S_i$ is an
 $S_{++}$-filter-regular element with respect to  $M,$ we have
$$\ell_A\Big[\Big(\dfrac{M}{xM}\Big)_{(n_1,\ldots,n_d)}\Big]=
\ell_A[M_{(n_1,\ldots,n_{d})}]-\ell_A[M_{(n_1,\ldots,n_i-1,\ldots,n_{d})}]\eqno(7)$$ for all large $n_1,\ldots,n_d$ by \cite[Remark 2.6]{VM}.
Denote by $P(n_1,\ldots,n_i,\ldots,n_d)$ the polynomial of $\ell_A[M_{(n_1,\ldots,n_{d})}]$ and
$Q(n_1,\ldots,n_d)$ the polynomial of $\ell_A\Big[\Big(\dfrac{M}{xM}\Big)_{(n_1,\ldots,n_d)}\Big],$ from (7) we have $$Q(n_1,\ldots,n_d)=P(n_1,\ldots,n_i,\ldots,n_d)-
P(n_1,\ldots,n_i-1,\ldots,n_d).\eqno(8)$$
Since  $e(M;k_1,\ldots,k_d)\ne 0 $ and $k_i > 0,$  by (8)
 we get $\deg Q =\deg P -1$ and  $$e(M;h_1,\ldots,h_d)=e\Big(\dfrac{M}{xM};h_1,\ldots,h_i-1,\ldots,h_d\Big)\; \text {for all } h_i > 0. \eqno(9) $$
 By $(9)$,
$$\begin{array}{l}\sum_{\mid\mathrm{\bf h}\mid\:=\:\ell-1;\; h_i >0}e(M;\mathrm{\bf h})
=\sum_{\mid \mathrm{\bf h}\mid\:=\:\ell-1;\; h_i >0}e\Big(\dfrac{M}{xM};h_1,\ldots,h_i-1,\ldots, h_d\Big).\end{array}$$
Since
$\widetilde{e}\Big(\dfrac{M}{xM}\Big) = \sum_{\mid \mathrm{\bf h}\mid\:=\:\ell-1;\; h_i >0}e\Big(\dfrac{M}{xM};h_1,\ldots,h_i-1,\ldots, h_d\Big),$
$$\widetilde{e}\Big(\dfrac{M}{xM}\Big) = \sum_{\mid \mathrm{\bf h}\mid\:=\:\ell-1;\; h_i >0}e(M;\mathrm{\bf h}).$$
We have (i). Remember that
$$\begin{array}{l}\widetilde{e}(M)=\sum_{\mid \mathrm{\bf h}\mid\:=\:\ell-1}
e(M;\mathrm{\bf h})
=\sum_{\mid\mathrm{\bf h}\mid\:=\:\ell-1;\; h_i >0}e(M;\mathrm{\bf h})
+\sum_{\mid \mathrm{\bf h}\mid\:=\:\ell-1;\; h_i =0}e(M;\mathrm{\bf h}).
\end{array}$$
Thus,
$$\widetilde{e}(M) = \widetilde{e}\Big(\dfrac{M}{xM}\Big)+\sum_{\mid\mathrm{\bf h}\mid\:=\:\ell-1;\; h_i =0}e(M;\mathrm{\bf h}).\eqno(10)$$
Now, we prove (ii).
 Choose $v \gg 0$ such that $$P(n_1,\ldots,n_d)= \ell_A[M_{(n_1,\ldots,n_{d})}]$$ for all $n_1,\ldots,n_d \ge v.$ Then $P(n_1,\ldots,v, \ldots,n_d)= \ell_A[M_{(n_1,\ldots,v,\ldots,n_{d})}]$ for all $$n_1,\ldots,n_{i-1},n_{i+1},\ldots, n_d \ge v\; \text{and}\; n_i = v.$$
Note that $$\ell_A[M_{(n_1,\ldots,v,\ldots,n_{d})}] = \ell_A[S_i^v{M_{\widehat{i}}}_{(n_1,\ldots,0,\ldots,n_{d})}]$$ and the terms of total degree $\ell-1$ in the polynomial $$P(n_1,\ldots,v, \ldots,n_d)= \ell_A[S_i^v{M_{\widehat{i}}}_{(n_1,\ldots,0,\ldots,n_{d})}]$$ have the form
$$ \sum_{h_1\:+\:\cdots+0+\:\cdots+\:h_d\;=\;\ell-1}e(M;h_1,\ldots,0,\ldots, h_d)\dfrac{n_1^{h_1}\cdots v^0\cdots n_d^{h_d}}{h_1!\cdots 0!\cdots h_d!}.$$
This follows that
$\sum_{\mid\mathrm{\bf h}\mid\:=\:\ell-1;\; h_i =0}e(M;\mathrm{\bf h})\ne 0$ if and only if
$\dim_{{S_{\widehat{i}}}^\triangle}[S_i^vM_{\widehat{i}}]^\triangle =\ell$ for some $v \gg 0.$ In this case,
$$e(M;h_1,\ldots,h_{i-1},0,h_{i+1},\ldots, h_d) = e(S_i^vM_{\widehat{i}};h_1,\ldots,h_{i-1},h_{i+1},\ldots, h_d)$$ for all
$v \gg 0$ by Remark 5.1.  Therefore
$\widetilde{e}\big({S_i^vM_{\widehat{i}}}\big)=
\sum_{\mid\mathrm{\bf h}\mid\:=\:\ell-1;\; h_i =0}e(M;\mathrm{\bf h})$ for all  $v \gg 0.$  (ii) is proved.
By $(10)$ and (ii) we immediately get (iii) and (iv). $\blacksquare$
\end{proof}
\vskip 0.2cm
\vspace*{10pt}
 We now will discuss  how  particular cases of Theorem  5.2 can be treated.\\
\noindent
Remember that if the multiplicities of $M;$ $\dfrac{M}{xM};$ $S_i^vM_{\widehat{i}}$ are expressed as the sums of all the  mixed multiplicities of $M;$ $\dfrac{M}{xM};$ $S_i^vM_{\widehat{i}},$ respectively   then
$$e(M) = \widetilde{e}(M);\; e\Big(\dfrac{M}{xM}\Big)=\widetilde{e}\Big(\dfrac{M}{xM}\Big);\; e(S_i^vM_{\widehat{i}})= \widetilde{e}(S_i^vM_{\widehat{i}}).$$\\
 Hence
  as  an immediate consequence of Theorem 5.2, we have the following result.
\vskip 0.2cm
 \noindent
{\bf Corollary 5.3.}\;
{\it  Set $\dim_{S^\triangle} M^\triangle = \ell.$  Assume that  $\ell > 1.$ Then we have:
\begin{itemize}
 \item[$\mathrm{(i)}$] If $e(M)$ is the sum of all the  mixed multiplicities of $M$ and $\dim_{{S_{\widehat{i}}}^\triangle}[S_i^vM_{\widehat{i}}]^\triangle =\ell$ for some $v \gg 0$ then   $e(M)= \widetilde{e}\Big(\dfrac{M}{xM}\Big)+\widetilde{e}
({S_i^vM_{\widehat{i}}}) $ $\text{ for all}\;\;  v \gg 0.$
 \item[$\mathrm{(ii)}$] If $e\Big(\dfrac{M}{xM}\Big)$ is the sum of all the   mixed multiplicities of $\dfrac{M}{xM}$  then
$$e\Big(\dfrac{M}{xM}\Big) = \sum_{\mid\mathrm{\bf h}\mid\:=\:\ell-1;\; h_i >0}e(M;\mathrm{\bf h}).$$
\item[$\mathrm{(iii)}$] If $\dim_{{S_{\widehat{i}}}^\triangle}[S_i^vM_{\widehat{i}}]^\triangle =\ell$ and $e(S_i^vM_{\widehat{i}})$ is the sum of all the  mixed multiplicities of $S_i^vM_{\widehat{i}}$ for some $v \gg 0$ then
     $$ e(S_i^vM_{\widehat{i}})=
\sum_{\mid\mathrm{\bf h}\mid\:=\:\ell-1;\; h_i =0}e(M;\mathrm{\bf h})
 \;\text {for all}\; v \;\gg 0.$$
\end{itemize}}

Remember that  $\mathbb{S}= F_J(J,\mathrm{\bf I}; R)$ and
$\mathbb{M}= F_J(J,\mathrm{\bf I};N);$ $ I = I_1\cdots I_d$  is not contained in $ \sqrt{\mathrm{Ann}{N}};$ $\dim \dfrac{N}{0_N:{ I}^\infty} = q.$  For any  $i=1,\ldots,d,$ set $$\mathfrak R(\mathrm{\bf I}_{\widehat{i}}\,; N) = \mathfrak R(I_1,\ldots,I_{i-1},I_{i+1},\ldots,I_d;N).$$
    Recall that   by Proposition  4.5(iii), $$\mathbb{M}_{\widehat{i}}
\cong F_J(J, I_1,\ldots,I_{i-1},I_{i+1},\ldots,I_d;N).$$
Upon simple computation, we get $$ \mathbb{S}_i^v\mathbb{M}_{\widehat{i}}
\cong F_J(J, I_1,\ldots,I_{i-1},I_{i+1},\ldots,I_d;I_i^vN).$$
Set
$$\overline{N} = \dfrac{N}{0_N: {I}^\infty}\;\; \text{ and}\;
\overline{\mathbb{M}}= F_J(J,\mathrm{\bf I};\overline{N}).$$
Then since  $\text{ht} \dfrac{ I+\text{Ann}\overline{N}}{\text{Ann}\overline{N}} > 0,$  we have
$$ \dim I_i^v\overline{N} = \dim \overline{N} > \dim \dfrac{\overline{N}}{I_i^v\overline{N}}$$ for any $1\leqslant i\leqslant d$ and for all $v >0.$
Hence from short exact sequences
$$ 0\longrightarrow I_i^v\overline{N} \longrightarrow \overline{N}\longrightarrow \dfrac{\overline{N}}{I_i^v\overline{N}}\longrightarrow 0, $$ by Corollary 3.9(ii)(b)
we get $$e\big(\big(J,\mathfrak R(\mathrm{\bf I}_{\widehat{i}}\,; R)_+\big) ; \mathfrak R(\mathrm{\bf I}_{\widehat{i}}\,; \overline{N})\big) = e\big(\big(J, \mathfrak R(\mathrm{\bf I}_{\widehat{i}}\,; R)_+\big) ;\mathfrak R(\mathrm{\bf I}_{\widehat{i}}\,; I_i^v\overline{N})\big).$$ On the other hand
 $$\widetilde{e}
\big(\mathbb{S}_i^v{\mathbb{M}}_{\widehat{i}}\big)=  e\big(\big(J,\mathfrak R(\mathrm{\bf I}_{\widehat{i}}\,; R)_+\big); \mathfrak R(\mathrm{\bf I}_{\widehat{i}}\,; I_i^v\overline{N})\big)$$ by Corollary 2.6. Hence
$$\widetilde{e}
\big(\mathbb{S}_i^v{\mathbb{M}}_{\widehat{i}}\big)=  e\big(\big(J,\mathfrak R(\mathrm{\bf I}_{\widehat{i}}\,; R)_+\big) ;\mathfrak R(\mathrm{\bf I}_{\widehat{i}}\,; \overline{N})\big).$$
This fact yields:
\vskip 0.2cm
\noindent {\bf Note 5.4.} We have
$$\widetilde{e}
\big(\mathbb{S}_i^v{\mathbb{M}}_{\widehat{i}}\big)=  e\big(\big(J,\mathfrak R(\mathrm{\bf I}_{\widehat{i}}\,; R)_+\big) ; \mathfrak R(\mathrm{\bf I}_{\widehat{i}}\,;\overline{N})\big).$$

\vskip 0.2cm
Put $\frak J = (J,\mathfrak R(\mathrm{\bf I};R)_+)$ and
$\frak J_{\widehat{i}}= (J,\mathfrak R(\mathrm{\bf I}_{\widehat{i}}\,; R)_+).$
Then  as a  consequence of Theorem 5.2 and Proposition  4.5  we obtain the following results.
\vskip 0.2cm
\noindent {\bf Theorem 5.5.}\; {\it Assume that $e(J^{[k_0+1]},\mathrm{\bf I}^{[\mathrm{\bf k}]};N) \ne 0$ and $k_i > 0.$
Let $x \in I_i$ be a weak-$(FC)$-element of $N$ with respect to $(J, I_1,\ldots, I_d).$ Then}

\begin{itemize}
\item[$\mathrm{(i)}$]\; $e\Big({\frak J};\mathfrak R\big(\mathrm{\bf I};\dfrac{N}{xN:{ I}^\infty}\big)\Big)=
\sum_{h_0+\mid\mathrm{\bf h}\mid =\:q-1;\; h_i >0}e\big(J^{[h_0+1]},\mathrm{\bf I}^{[\mathrm{\bf h}]}; N\big).$
\item[$\mathrm{(ii)}$] \; $e\Big(\frak J; \mathfrak R\big(\mathrm{\bf I};\dfrac{N}{0_N:{ I}^\infty}\big)\Big)= e\Big(\frak J; \mathfrak R\big(\mathrm{\bf I};\dfrac{N}{xN:{ I}^\infty}\big)\Big)+ e\Big(\frak J_{\widehat{i}}; \mathfrak R\big(\mathrm{\bf I}_{\widehat{i}}\,;\dfrac{N}{0_N:{ I}^\infty}\big)\Big).$
\item[$\mathrm{(iii)}$]\; $e\Big(\mathfrak R\big(\mathrm{\bf I};\dfrac{N}{0_N:{ I}^\infty}\big)\Big)= e\Big(\mathfrak R\big(\mathrm{\bf I};\dfrac{N}{xN:{ I}^\infty}\big)\Big)+
e\Big(\mathfrak R\big(\mathrm{\bf I}_{\widehat{i}}\,;\dfrac{N}{0_N:{I}^\infty}\big)\Big).
$
\end{itemize}

\begin{proof}\; Denote by
$\bar x$ the image of $x$ in $\mathbb{S}_i.$
 Since $x \in I_i$ is a weak-(FC)-element of $N$ with respect to $(J, I_1,\ldots, I_d),$ $\bar x$ is an $\mathbb{S}_{++}$-filter-regular element with respect to $\mathbb{M}$ by Proposition  4.5(i). By Proposition  4.5(ii),
$\widetilde{e}\big(\mathbb{M}/\bar x\mathbb{M}\big) =\widetilde{e}\Big(F\big(J,\mathrm{\bf I};\dfrac{{N}}{x{N}}\big)\Big).$
Hence
 \begin{align*}
\widetilde{e}\big(\mathbb{M}/\bar x\mathbb{M}\big) &=e\Big(F_J\big(J,\mathrm{\bf I};\dfrac{N}{xN: I^\infty}\big)\Big)\\&= e\Big(\big(J,\mathfrak R(\mathrm{\bf I}; R)_+\big); \mathfrak R\big(\mathrm{\bf I}; \dfrac{N}{xN: I^\infty}\big)\Big)\\&= e\Big({\frak J};\mathfrak R\big(\mathrm{\bf I};\dfrac{N}{xN:{ I}^\infty}\big)\Big)  \end{align*} by   Corollary 2.6.
Thus, we get(i) by Theorem 5.2(i).
 Now, since
$$ \widetilde{e}\big(\mathbb{M}\big)=  e\Big(\frak J; \mathfrak R\big(\mathrm{\bf I};\dfrac{N}{0_N:{ I}^\infty}\big)\Big)\;\; \text{and}\;\;
\widetilde{e}\big(\mathbb{M}/\bar x\mathbb{M}\big) =e\Big(\frak J; \mathfrak R\big(\mathrm{\bf I};\dfrac{N}{xN:{ I}^\infty}\big)\Big)$$
by Corollary 2.6, and
$$\widetilde{e}
\big(\mathbb{S}_i^v{\mathbb{M}}_{\widehat{i}}\big)=  e\big(\big(J,\mathfrak R(\mathrm{\bf I}_{\widehat{i}}\,; R)_+\big) ;\mathfrak R(\mathrm{\bf I}_{\widehat{i}}\,; \overline{N})\big)= e\Big(\frak J_{\widehat{i}}; \mathfrak R(\mathrm{\bf I}_{\widehat{i}}\,; \dfrac{N}{0_N:{ I}^\infty})\Big)$$ by Note 5.4,
 we have (ii) by Theorem 5.2(iii). Choose $J = \frak n,$ we get (iii) by (ii).
 $\blacksquare$
\end{proof}

   Remember that if ht $\dfrac{I+\text{Ann}N}{\text{Ann}N} > 0,$ then
$$e\Big({\frak J};\mathfrak R\big(\mathrm{\bf I};\dfrac{N}{0_N:{ I}^\infty}\big)\Big)= e\big({\frak J};\mathfrak R\big(\mathrm{\bf I}; N\big)\big)$$
by Remark 2.7. Hence as an immediate consequence of Theorem 5.5, we have the following result.

\noindent {\bf Corollary 5.6.}\; {\it Assume that $\mathrm{ht}$ $\dfrac{I+\mathrm{Ann}N}{\mathrm{Ann}N} > 0;$ $e(J^{[k_0+1]},\mathrm{\bf I}^{[\mathrm{\bf k}]};N) \ne 0$ and $k_i > 0.$
Let $x \in I_i$ be a weak-$(FC)$-element of $N$ with respect to $(J, I_1,\ldots, I_d).$ Then}

\begin{itemize}
\item[$\mathrm{(i)}$]\; $e\Big(\frak J; \mathfrak R\big(\mathrm{\bf I};\dfrac{N}{xN:{ I}^\infty}\big)\Big)=
\sum_{h_0+\mid\mathrm{\bf h}\mid =\:q-1;\; h_i >0}e\big(J^{[h_0+1]},\mathrm{\bf I}^{[\mathrm{\bf h}]}; N\big).$
\item[$\mathrm{(ii)}$] \; $e\big(\frak J; \mathfrak R\big(\mathrm{\bf I}; N\big)\big)= e\Big(\frak J; \mathfrak R\big(\mathrm{\bf I};\dfrac{N}{xN:{ I}^\infty}\big)\Big)+ e\big(\frak J_{\widehat{i}} ;  \mathfrak R\big(\mathrm{\bf I}_{\widehat{i}}\,; {N}\big)\big).$
\item[$\mathrm{(iii)}$]\; $e\big(\mathfrak R\big(\mathrm{\bf I}; {N}\big)\big)= e\Big(\mathfrak R\big(\mathrm{\bf I};\dfrac{N}{xN:{ I}^\infty}\big)\Big)+
e\big(\mathfrak R\big(\mathrm{\bf I}_{\widehat{i}}\,; {N}\big)\big).
$
\end{itemize}

Suppose that $e\big(J^{[k_0+1]},\mathrm{\bf I}^{[\mathrm{\bf k}]};N\big) \ne 0$  and
$x_1,\ldots, x_p$ $(p \leqslant k_i)$ is  a  weak-(FC)-sequence in $I_i.$
By Theorem 5.5 and by induction on $p$, we get the following corollary.
\vskip 0.2cm
\noindent {\bf Corollary 5.7.}\;{\it Let $e(J^{[k_0+1]},\mathrm{\bf I}^{[\mathrm{\bf k}]};N) \ne 0$ and $x_1,\ldots, x_p \in I_i$ $(p \leqslant k_i)$ be a  weak-$(FC)$-sequence  of $N$ with respect to $(J, I_1,\ldots, I_d).$  Then}
 \begin{align*}
 e\Big(\frak J; \mathfrak R\big(\mathrm{\bf I};\dfrac{N}{0_N: I^\infty}\big)\Big)
 &= e\Big(\frak J; \mathfrak R\big(\mathrm{\bf I};\dfrac{N}{(x_1,\ldots, x_p)N:{ I}^\infty}\big)\Big)\\ &+ \sum_{j=0}^{p-1}e\Big(\frak J_{\widehat{i}}; \mathfrak R\big(\mathrm{\bf I}_{\widehat{i}}\,;\dfrac{N}{(x_1,\ldots, x_j)N:{I}^\infty}\big)\Big). \end{align*}

In particular, if $d=1$ then $I= I_1.$ Put $p = \max\{i \;|\; e(J^{[q-i]},I^{[i]};N) \ne 0\}$ and
assume that  $x_1,\ldots, x_p$ is a  weak-(FC)-sequence of $N$ with respect to $(J,I).$
Then   by \cite{Vi1, Vi2}(see \cite[Proposition 3.3(iii) and Theorem 3.4(iii)]{MV}), $x_1,\ldots, x_p$ is a  maximal (FC)-sequence of $N$ with respect to $(J,I).$   By Corollary 5.7,     $$e\Big(\frak J; \mathfrak R( I;\dfrac{N}{0:{ I}^\infty})\Big)= e\Big(\frak J; \mathfrak R(I;\dfrac{N}{(x_1,\ldots, x_p)N:{I}^\infty})\Big)+ \sum_{i=0}^{p-1}e\Big( J;\dfrac{N}{(x_1,\ldots,x_i)N:{ I}^\infty}\Big).$$  Since \;$p$ \;is  maximal, \; $e(J^{[q-i]},I^{[i]};N) \ne 0$ if \;and \;only if $ 0 \leqslant i \leqslant p$ \;by \cite{Vi1}. Consequently by \cite{Vi1}(see \cite[Proposition 3.3(i)]{MV}),
$$e(J^{[q-p-i]},I^{[p+i]};N)= e\Big(J^{[q-p-i]},I^{[i]};\dfrac{N}{(x_1,\ldots, x_p)N}\Big) \ne 0$$ if and only if $i=0$.
Therefore  by Corollary 2.5(ii),

$$ e\Big(\frak J; \mathfrak R(I;\dfrac{N}{(x_1,\ldots, x_p)N:{I}^\infty})\Big)= e\Big(J^{[q-p]},I^{[0]};\dfrac{N}{(x_1,\ldots, x_p)N}\Big).$$  On \;the other \; hand \;
$e\Big(J^{[q-p]},I^{[0]};\dfrac{N}{(x_1,\ldots, x_p)N}\Big)= e\Big(J;\dfrac{N}{(x_1,\ldots, x_p)N:{I}^\infty}\Big)$ by \cite[Lemma 3.2]{MV}. \noindent Hence
$e\Big(\frak J;\mathfrak R(I;\dfrac{N}{(x_1,\ldots, x_p)N:{I}^\infty})\Big)=
e\Big(J;\dfrac{N}{(x_1,\ldots, x_p)N:{I}^\infty}\Big).$
Thus, $$e\Big(\frak J; \mathfrak R(I;\dfrac{N}{0_N: { I}^\infty})\Big)=\sum_{j=0}^pe\Big(J;\dfrac{N}{(x_1,\ldots, x_j)N:{I}^\infty}\Big).$$  Then  we have the following corollary.
\vskip 0.2cm
\noindent {\bf Corollary 5.8.}\;
 $e\Big(\frak J; \mathfrak R(I;\dfrac{N}{0_N: { I}^\infty})\Big)=\sum_{j=0}^pe\Big(J;\dfrac{N}{(x_1,\ldots, x_j)N:{I}^\infty}\Big).$\\

In the case that $$\text{ht}\dfrac{ I+\text{Ann}N}{\text{Ann}N}>0,\;
e\Big(\frak J; \mathfrak R(I;\dfrac{N}{0_N: { I}^\infty})\Big)= e\big(\frak J; \mathfrak R(I;N)\big)$$ by Remark 2.7.  We get the following result which  is proved by \cite{MV}.
\vskip 0.2cm
\noindent {\bf Corollary 5.9} \cite[Theorem 4.2]{MV}.\; {\it If\;  $\mathrm{ht}\dfrac{ I+\mathrm{Ann}N}{\mathrm{Ann}N}>0,$    then
 $$e\big(\frak J; \mathfrak R(I;N)\big)=\sum_{j=0}^pe\Big(J;\dfrac{N}{(x_1,\ldots, x_j)N:{I}^\infty}\Big).$$}

\end{document}